\documentclass[11pt]{article}

\usepackage[T1]{fontenc}
\usepackage{lmodern}
\usepackage[margin=1in]{geometry}
\usepackage{amsmath,amssymb,amsthm,aliascnt}
\usepackage{booktabs,braket,xcolor}
\usepackage[numbers,sort&compress]{natbib}
\usepackage[hidelinks]{hyperref}
\usepackage[capitalize,noabbrev]{cleveref}

\numberwithin{equation}{section}
\theoremstyle{plain}
\newtheorem{theorem}{Theorem}[section]
\newaliascnt{lemma}{theorem}
\newtheorem{lemma}[lemma]{Lemma}
\aliascntresetthe{lemma}
\newaliascnt{corollary}{theorem}
\newtheorem{corollary}[corollary]{Corollary}
\aliascntresetthe{corollary}
\theoremstyle{remark}
\newaliascnt{remark}{theorem}
\newtheorem{remark}[remark]{Remark}
\aliascntresetthe{remark}
\crefname{theorem}{Theorem}{Theorems}
\crefname{lemma}{Lemma}{Lemmas}
\crefname{corollary}{Corollary}{Corollaries}
\crefname{remark}{Remark}{Remarks}

\newenvironment{keywords}
  {\par\smallskip\noindent\textbf{Keywords.}\ }
  {\par\smallskip}
\newenvironment{MSCcodes}
  {\par\smallskip\noindent\textbf{Mathematics Subject Classification.}\ }
  {\par\smallskip}

\hypersetup{
  pdftitle={Finite Element Approximation of Three-Dimensional Elliptic Optimal Control Problems with Curve-Supported Controls},
  pdfauthor={D. Leykekhman and B. Vexler},
  pdfsubject={Elliptic optimal control with curve-supported controls},
  pdfkeywords={elliptic optimal control, curve-supported controls, finite element approximation, variational discretization, weighted error estimates}
}

\newcommand{\abs}[1]{\lvert#1\rvert}

\newcommand{\half}{\frac{1}{2}}

\newcommand{\na}{\nabla}

\newcommand{\Om}{\Omega}
\newcommand{\R}{\mathbb{R}}

\newcommand{\Hg}{\mathcal{H}^{1}}

\newcommand{\vertiii}[1]{{\left\vert\kern-0.25ex\left\vert\kern-0.25ex\left\vert #1
    \right\vert\kern-0.25ex\right\vert\kern-0.25ex\right\vert}}

\definecolor{darkred}{rgb}{.7,0,0}

\definecolor{green}{rgb}{0,0.7,0}

\title{Finite Element Approximation of Three-Dimensional Elliptic
Optimal Control Problems with Curve-Supported Controls}
\author{D. Leykekhman \and B. Vexler}
\date{}

\begin{document}
\maketitle

\begin{abstract}
We consider a linear-quadratic elliptic optimal control problem in three dimensions with a box-constrained control supported on a compact embedded $C^{1,1}$
 curve strictly contained in a bounded convex polyhedral domain. The right-hand side of the state equation is a measure supported on the curve, and the state is singular along it. We establish well-posedness and first-order optimality conditions and derive curve-adapted weighted inequalities and finite element error estimates for the singular state. For continuous piecewise linear finite elements on quasi-uniform tetrahedral meshes and variational discretization of the control,
assuming exact integration over the curve, we prove the almost quadratic estimate
$\|\bar q-\bar q_h\|_{L^2(\gamma)}
\le Ch^2(1+\abs{\ln h})^3$.
The proof of this estimate requires new results, which are of independent interest. The convergence rate observed in a numerical experiment is consistent with the predicted algebraic order.
\end{abstract}




\section{Introduction}\label{sec: intro}
We consider the following optimal control problem
\begin{subequations}
\begin{equation}\label{eq: minimize}
\text{Minimize}\quad J(q,u)=\frac{1}{2}\|u-u_d\|^2_{L^2(\Omega)}+\frac{\alpha}{2}\|q\|^2_{L^2(\gamma)},
\end{equation}
subject to
\begin{equation}\label{eq: the problem}
\begin{aligned}
-\Delta u  &=q\delta_\gamma \quad \text{in $\Omega$}\\
u &=0  \quad \text{on $\partial\Omega$}
\end{aligned}
\end{equation}
and
\begin{equation}\label{eq: the problem_control_constraints}
q_a\le q\le q_b\ \text{a.e. on }\gamma
\end{equation}
\end{subequations}
Here $\Omega\subset\mathbb R^3$ is a bounded convex polyhedral domain,
and $\gamma\subset\Om_0\Subset\Omega$ is a compact embedded $C^{1,1}$
curve with two endpoints. The regularization parameter satisfies
$\alpha>0$, the control bounds $q_a,q_b\in\mathbb R$ satisfy $q_a<q_b$,
and the desired state satisfies $u_d\in L^\infty(\Omega)$.

We denote by $\delta_\gamma:=\mathcal H^1\llcorner\gamma$ the arclength
measure on $\gamma$. The control space is $Q:=L^2(\gamma,\Hg)$,
and the set of admissible controls is
\[
Q_{\rm ad}:=\{q\in Q:q_a\le q\le q_b\ \text{a.e. on }\gamma\}.
\]
The right-hand side of the state equation is the measure
$q\delta_\gamma$, defined by
\[
\langle q\delta_\gamma, \varphi \rangle = \int_\gamma q\varphi\ d\Hg
\]
for $\varphi$ with appropriate regularity, see below.
  
Problems involving line sources have attracted considerable
interest because of their applications to 3D--1D coupling
\cite{BERRONE_2023}.

Although $q\delta_\gamma\notin H^{-1}(\Omega)$ in general, the
state equation is well posed in $W_0^{1,p}(\Omega)$ for
$6/5<p<2$, as shown in Section~\ref{sec: existence}. Its conforming finite element
discretization is nevertheless well defined.
 Finite element approximations of elliptic problems with line
 sources were studied in weighted Sobolev spaces in
 \cite{DAngeloC_2012, DrelichmanI_DuranR_OjeaI_2020}.
 A singularity splitting approach combined with finite elements
 was developed in
 \cite{GjerdeI_KumarK_NordbottenJ_WohlmuthB_2019},
 while a singularity splitting deep Ritz method was proposed in
 \cite{HuT_JinB_ZhouZ_2023}. A discontinuous Galerkin discretization of elliptic and parabolic
problems with line sources is analyzed in \cite{MasriR_ShenB_RiviereB_2023}.
More recently, Gao and Huang \cite{gao2026optimallocalerrorestimates}
studied Lagrange finite element approximations of elliptic problems
with measure data. Using interior estimates, they obtained 
$L^2$- and $H^1$-error bounds on subdomains at positive distance
from the support of the source.

A closely related class of elliptic optimal control problems with
controls acting on lower-dimensional manifolds was analyzed in
\cite{GongW_WangG_YanN_2014}. For a curve-supported control in three
dimensions, that analysis yields a control-error estimate of almost
order $h^{3/2}$.

We use continuous piecewise linear finite elements on quasi-uniform,
shape-regular tetrahedral meshes and variational discretization of the
control. Assuming exact integration over the curve, we prove
\[
\|\bar q-\bar q_h\|_{L^2(\gamma)}
\le Ch^2(1+\abs{\ln h})^3,
\]
where $\bar q$ and $\bar q_h$ denote the optimal controls of the
continuous and discrete problems, respectively.

The almost quadratic convergence rate is analogous to that
for point-supported controls in two dimensions;
see \cite[Chapter~10]{vexler2025numerical}.
In both settings the control acts on a set of codimension two.
Our analysis relies on several estimates that are also of independent interest for elliptic equations with curve-supported sources. Writing
\[
d_\gamma(x):=\operatorname{dist}(x,\gamma),
\qquad
\sigma_h(x):=\sqrt{d_\gamma(x)^2+h^2},
\qquad
\ell_h:=1+|\ln h|,
\]
we first establish the weighted Hardy inequality
\[
\|\sigma_h^{-1}v\|_{L^2(\Omega)}
\le C\ell_h\|\nabla v\|_{L^2(\Omega)}
\qquad\text{for all }v\in H^1_0(\Omega);
\]
see \cref{lem:weighted-hardy}. We further prove the discrete Sobolev-type inequality
\[
\|v_h\|_{L^2(\gamma)}
\le C\ell_h^{1/2}\|\nabla v_h\|_{L^2(\Omega)}
\qquad\text{for all }v_h\in V_h;
\]
see \cref{lem:discrete-sobolev}.
Although restriction to $\gamma$ is not a bounded operator
from $H^1(\Omega)$ to $L^2(\gamma)$, this estimate provides
control of the trace for finite element functions, with
only a logarithmic dependence on the mesh size.
It is analogous to the classical discrete Sobolev
inequality in two dimensions
\cite[Lemma~4.9.2]{BrennerScott2008}.

For the state equation with an arbitrary density
$q\in L^2(\gamma)$, let $u$ and $u_h$ denote its continuous
and finite element solutions, respectively.
We prove the global error estimate
\[
\|u-u_h\|_{L^2(\Omega)}
\le Ch\|q\|_{L^2(\gamma)};
\]
see \cref{lem:state-L2-error}. For straight line sources, global $L^2$-error estimates of order
$h^{1-\varepsilon}$, for arbitrarily small $\varepsilon>0$
and with constants depending on $\varepsilon$, were established
for continuous Galerkin approximations in
\cite[Lemma~3.3]{GongW_WangG_YanN_2014}
and for discontinuous Galerkin approximations in
\cite[Theorem~1]{MasriR_ShenB_RiviereB_2023}.
Our estimate removes this loss in the convergence order.

A dyadic decomposition into tubular layers, combined with local pointwise finite element estimates and integration along $\gamma$, yields the weighted error estimate
\[
\|\sigma_h(u-u_h)\|_{L^2(\Omega)}
\le Ch^2\ell_h^{3/2}\|q\|_{L^2(\gamma)};
\]
see \cref{thm: weighted error in L2}.
As consequences, we obtain
\[
\|u-u_h\|_{L^1(\Omega)}
\le Ch^2\ell_h^2\|q\|_{L^2(\gamma)}
\]
and
\[
\|u-u_h\|_{L^2(\Omega_d)}
\le Ch^2\ell_h^{3/2}d^{-1}\|q\|_{L^2(\gamma)},
\qquad
\Omega_d:=\{x\in\Omega:d_\gamma(x)>d\};
\]
see \cref{cor: L1}.
Thus, the approximation converges almost quadratically
both in the global $L^1$ norm and in the $L^2$ norm away
from the curve. The global $L^1$ estimate recovers the rate established
for straight line sources in
\cite[Lemma~3.6]{GongW_WangG_YanN_2014}. The latter estimate holds with a constant
independent of $h$ and $d>0$, making the deterioration
as the observation region approaches the curve explicit.
It complements local estimates on regions at a fixed
positive distance from the source, such as
\cite[Theorem~4.1]{gao2026optimallocalerrorestimates}.

Finally, combining the weighted state error estimate with
the weighted Hardy and discrete Sobolev inequalities,
and with local pointwise estimates for the adjoint,
yields the stated control error estimate.

The regularity results and error estimates extend, with unchanged
convergence orders, to simple curves consisting of finitely many
regular $C^{1,1}$ arcs, including non-self-intersecting polygonal
chains, provided that the curve remains at positive distance
from $\partial\Omega$.
This follows by applying the geometric and analytic estimates
on the individual arcs and combining their contributions.

The rest of the paper is organized as follows.
Section~2 collects geometric properties of the control curve.
Section~3 addresses well-posedness, optimality conditions,
and regularity.
Section~4 establishes auxiliary inequalities and finite element
error estimates for the state equation.
The control error analysis is presented in Section~5,
followed by a numerical experiment in Section~6.
The proofs of the tubular coordinate and integration lemmas
are collected in the appendix.
\section{Preliminaries on the curve}
Let $\gamma\subset\Omega_0\Subset\Omega$ be a compact embedded
$C^{1,1}$ curve with two endpoints. We fix an injective arclength
parametrization $\Gamma\in C^{1,1}([0,L];\mathbb R^3)$, so that
\begin{equation}\label{eq:parametrization_curve}
\gamma=\Gamma([0,L]), \qquad \abs{\Gamma'(s)}=1
\quad\text{for all }s\in[0,L],
\end{equation}
where $L=\mathcal H^1(\gamma)$. We write
\begin{equation}\label{eq:dist_and_tab_neigbourhood}
d_\gamma(x):=\operatorname{dist}(x,\gamma), \qquad
T_r(\gamma):=\{x\in\mathbb R^3:d_\gamma(x)<r\}.
\end{equation}
Since $\gamma$ is compactly contained in $\Omega$, we have
$d_{\gamma,\partial\Omega}:=\operatorname{dist}(\gamma,\partial\Omega)>0$.

\begin{lemma}[Local length estimates]\label{lem:curve-local-length}
	There exist $r_{\rm len}>0$ and $C>0$, depending only on $\gamma$,
	such that for every $0<r\le r_{\rm len}$,
	\begin{align}
		\mathcal H^1(\gamma\cap B_r(y)) &\le Cr
		&&\text{for all }y\in\mathbb R^3, \label{eq:curve-length-upper}\\
		\mathcal H^1(\gamma\cap B_r(y)) &\ge r
		&&\text{for all }y\in\gamma. \label{eq:curve-length-lower}
	\end{align}
\end{lemma}
\begin{proof}
	By uniform continuity of $\Gamma'$, choose $0<\delta<L$ such that
	\[
	\abs{\Gamma'(s)-\Gamma'(t)}\le \half
	\]
	for all $0\le s<t\le L$ with $\abs{s-t}\le\delta$.
	For such $s$ and $t$, it follows by Cauchy-Schwarz inequality
	\[
	\begin{aligned}
	\abs{\Gamma(t)-\Gamma(s)}
	&\ge (\Gamma(t)-\Gamma(s))\cdot\Gamma'(s)\\
	&= \int_s^t \Gamma'(\xi)\cdot\Gamma'(s)\,d\xi\\
	&= \int_s^t \left(\Gamma'(s) +(\Gamma'(\xi)-\Gamma'(s))\right)\cdot\Gamma'(s)\,d\xi\\
	& \ge \int_s^t \left(1-\abs{\Gamma'(\xi)-\Gamma'(s)}\right)d\xi\\
	&\ge \tfrac12(t-s).
	\end{aligned}
	\]
	Moreover, compactness and injectivity give
	\[
	m_\delta:=
	\min_{\substack{s,t\in[0,L]\\ \abs{s-t}\ge\delta}}
	\abs{\Gamma(s)-\Gamma(t)}>0.
	\]
	Choose $0<r_{\rm len}\le L/2$ such that $2r_{\rm len}<m_\delta$.
	
	To prove \eqref{eq:curve-length-upper}, suppose that
	$\gamma\cap B_r(y)$ is nonempty and choose $s_0 \in [0,L]$ with 
	$\Gamma(s_0)\in\gamma\cap B_r(y)$. For every $s$ with
	$\Gamma(s)\in \gamma \cap B_r(y)$, we have
	\[
	\abs{\Gamma(s)-\Gamma(s_0)}<2r<m_\delta.
	\]
	Hence $\abs{s-s_0}<\delta$, and the local estimate above gives
	$\abs{s-s_0}\le 2\abs{\Gamma(s)-\Gamma(s_0)}<4r$. Thus $\gamma\cap B_r(y)$
	is contained in the image of a parameter interval of length at most
	$8r$. Since $\Gamma$ is parametrized by arclength, this proves
	\eqref{eq:curve-length-upper}.
	
	For \eqref{eq:curve-length-lower}, let $y=\Gamma(s_0)$. Since
	$r\le L/2$, at least one of the intervals $[s_0-r,s_0]$ and
	$[s_0,s_0+r]$ is contained in $[0,L]$. Its interior is mapped into
	$B_r(y)$ because $\abs{\Gamma(s)-\Gamma(s_0)}\le |s-s_0|$. The
	corresponding arc has length $r$.
\end{proof}

\begin{corollary}\label{cor:curve-ball-integral}
	For every $g\in L^2(\Omega)$ and $0<r\le r_{\rm len}$ with $r_{\rm len}$ from \cref{lem:curve-local-length}, we have
	\[
	\int_\gamma \|g\|_{L^2(B_r(x)\cap\Omega)}^2\,d\mathcal H^1(x) \le Cr\,\|g\|_{L^2(\Omega)}^2.
	\]
\end{corollary}

\begin{proof}
	By Fubini's theorem and the equivalence
	$y\in B_r(x)\iff x\in B_r(y)$, we obtain
	\[
	\begin{aligned}
		\int_\gamma \|g\|_{L^2(B_r(x)\cap\Omega)}^2\,d\mathcal H^1(x)
		&= \int_\gamma \int_\Omega |g(y)|^2
		\chi_{B_r(x)}(y)\,dy\,d\mathcal H^1(x)\\
		&= \int_\Omega |g(y)|^2
		\left(\int_\gamma \chi_{B_r(y)}(x)\,d\mathcal H^1(x)\right)dy\\
		&= \int_\Omega |g(y)|^2\,\mathcal H^1(\gamma\cap B_r(y))\,dy\\
		&\le Cr\,\|g\|_{L^2(\Omega)}^2,
	\end{aligned}
	\]
	where the last inequality follows from \eqref{eq:curve-length-upper}.
\end{proof}

\begin{lemma}\label{lem:curve-kernel-estimates}
	There exists a constant $C>0$ such that, for every
	$x\in\Omega\setminus\gamma$,
	\begin{align}
		\int_\gamma \frac{1}{|x-y|}\,d\mathcal H^1(y)
		&\le C\bigl(1+\abs{\ln d_\gamma(x)}\bigr),
		\label{eq1:curve-kernel-one}\\
		\int_\gamma \frac{1}{|x-y|^2}\,d\mathcal H^1(y)
		&\le \frac{C}{d_\gamma(x)}.
		\label{eq1:curve-kernel-two}
	\end{align}
	The constant $C$ depends only on $\gamma$ and $\Omega$.
\end{lemma}
\begin{proof}
	First, estimate \eqref{eq:curve-length-upper} extends to
	all $r>0$, after increasing the constant, since
	\[
	\mathcal H^1(\gamma\cap B_r(x))
	\le L
	\le \frac{L}{r_{\rm len}}r
	\qquad\text{for }r>r_{\rm len}.
	\]
	Fix $x\in\Omega\setminus\gamma$ and set
	$\rho=d_\gamma(x)>0$.
	Choose $R>\max\{1,\operatorname{diam}(\Omega)\}$ and let
	$N$ be the smallest positive integer such that $2^N\rho>R$.
	Then
	\[
	\gamma=\bigcup_{j=0}^{N-1} A_j,
	\qquad
	A_j:=\{y\in\gamma:2^j\rho\le |x-y|<2^{j+1}\rho\}.
	\]
	For $k\in\{1,2\}$, the upper length estimate gives
	\[
	\begin{aligned}
		\int_\gamma \frac{1}{|x-y|^k}\,d\mathcal H^1(y)
		&\le \sum_{j=0}^{N-1}
		(2^j\rho)^{-k}
		\mathcal H^1(\gamma\cap B_{2^{j+1}\rho}(x))\\
		&\le C\sum_{j=0}^{N-1}(2^j\rho)^{1-k}.
	\end{aligned}
	\]
	For $k=1$, the sum equals $N$, and
	\[
	N\le 1+\log_2(R/\rho)
	\le C(1+\abs{\ln\rho}).
	\]
	For $k=2$, we have
	\[
	\sum_{j=0}^{N-1}(2^j\rho)^{-1}
	\le \rho^{-1}\sum_{j=0}^{\infty}2^{-j}
	=\frac{2}{\rho}.
	\]
	This proves both estimates.
\end{proof}

For a nonempty closed set $A\subset\mathbb R^3$, its reach $\operatorname{reach}(A)$ is
defined as
\[
\sup\Set{r>0|\text{every }x\;\text{with}\;\operatorname{dist}(x,A)<r
\text{ has a unique nearest point in }A}.
\]
Since $\gamma$ is a compact embedded $C^{1,1}$ curve,
it has positive reach; see
\cite[Corollary~8.9(a)]{RatajZajicek2017}.

The next two lemmas provide tubular coordinates and an
integration estimate for neighborhoods of $\gamma$.
Their proofs, including the treatment of the endpoints,
are given in Appendix~\ref{app:curve}.

\begin{lemma}[Tubular coordinates around an extended curve]
	\label{lem:curve-tubular-coordinates}
	Let $\Gamma:[0,L]\to\gamma$ be the injective arclength
	parametrization fixed above.
	There exist $\delta>0$, an injective arclength extension
	\[
	\widetilde\Gamma\in
	C^{1,1}([-\delta,L+\delta];\mathbb R^3),
	\qquad
	\widetilde\Gamma|_{[0,L]}=\Gamma,
	\]
	with
	$\widetilde\gamma:=\widetilde\Gamma([-\delta,L+\delta])
	\Subset\Omega$, an open interval $J$ such that 
	\[
	[0,L]\subset J\Subset(-\delta,L+\delta),
	\]
	and Lipschitz continuous vector fields
	$n_1,n_2:\overline J\to\mathbb R^3$
	such that
	$(\widetilde\Gamma'(s),n_1(s),n_2(s))$
	is a positively oriented orthonormal basis of
	$\mathbb R^3$ for every $s\in\overline J$.
	
	For some $R>0$, the map
	\[
	\Psi(s,y)
	:=\widetilde\Gamma(s)+y_1n_1(s)+y_2n_2(s)
	\]
	is injective and Lipschitz continuous on
	$J\times B_R^{(2)}(0)$, where $B_R^{(2)}(0)$ denotes the open ball of radius $R$ centered at the origin in $\mathbb R^2$, and its image
	\[
	W^0:=\Psi(J\times B_R^{(2)}(0))
	\]
	is open and compactly contained in $\Omega$.
	Writing
	\[
	\widetilde d(x):=\operatorname{dist}(x,\widetilde\gamma),
	\]
	we have
	\begin{equation}\label{eq:tubular-distance}
		\widetilde d(\Psi(s,y))=|y|
		\qquad
		\text{for }s\in J,\quad |y|<R,
	\end{equation}
	and
	\begin{equation}\label{eq:tubular-distance-comparison}
		\widetilde d(x)\le d_\gamma(x)
		\qquad\text{for every }x\in\mathbb R^3.
	\end{equation}
	Moreover, there exist $r>0$ and an open set $W$ such that
	\begin{equation}\label{eq:tubular-cover}
		\overline{T_{2r}(\gamma)}
		\subset W\Subset W^0\Subset\Omega.
	\end{equation}
	
	The cylindrical coordinate map
	\[
	\Phi(s,\rho,\vartheta)
	:=\Psi(s,\rho\cos\vartheta,\rho\sin\vartheta)
	\]
	satisfies $|\partial_\rho\Phi|=1$ and
	\begin{equation}\label{eq:tubular-jacobian}
		\frac12\rho
		\le |\det D\Phi(s,\rho,\vartheta)|
		\le \frac32\rho
	\end{equation}
	for almost every
	$(s,\rho,\vartheta)\in J\times(0,R)\times(0,2\pi)$.
\end{lemma}

\begin{lemma}[Integration over tubular neighborhoods]
	\label{lem:curve-tubular-integration}
	There exist constants $r_0>0$ and $C>0$ such that
	$T_{r_0}(\gamma)\Subset\Omega$ and
	\begin{equation}\label{eq:curve-tubular-integration}
		\int_{T_r(\gamma)} F(d_\gamma(x))\,dx
		\le C\int_0^r F(\rho)\,\rho\,d\rho
	\end{equation}
	for every $0<r\le r_0$ and every nonnegative measurable
	function $F:[0,r_0]\to[0,\infty]$.
	In particular,
	\begin{equation}\label{eq:curve-tube-volume}
		|T_r(\gamma)|\le Cr^2
		\qquad\text{for }0<r\le r_0.
	\end{equation}
\end{lemma}

\section{Existence, optimality conditions, and regularity}\label{sec: existence}

First, for $6/5< p <2$, we  introduce a solution operator $S\colon Q=L^2(\gamma)\to W^{1,p}_0(\Om)$ which maps a given control $q\in Q$ to $u=u(q)=Sq\in W^{1,p}_0(\Om)$ by
\begin{equation}\label{eq: state}
(\na u,\na \varphi) = \langle q\delta_\gamma,\varphi\rangle :=\int_\gamma q\varphi\ d\Hg\quad \forall\varphi \in W^{1,p'}_0(\Om), 
\end{equation}
 where $\frac{1}{p}+\frac{1}{p'}=1$ and as a result  $2<p'<6$.

Since $\gamma$ is a compact embedded $C^{1,1}$ curve,
the trace theorem for submanifolds of codimension two
gives a bounded operator
\[
\operatorname{Tr}_\gamma:
W^{1,r}(\Omega)\longrightarrow L^r(\gamma),
\qquad 2<r<\infty;
\]
see, e.g., \cite[Chapter~VII, Theorem~1]{JonssonWallin1984}.
Since $\mathcal H^1(\gamma)<\infty$, H\"older's inequality
therefore yields
\[
\|\operatorname{Tr}_\gamma v\|_{L^2(\gamma)}
\le C_r\|v\|_{W^{1,r}(\Omega)}.
\]
In integrals over $\gamma$, Sobolev functions are understood
through their traces.
Taking $r=p'$, we obtain
\[
\left|\int_\gamma q\varphi\,d\mathcal H^1\right|
\le \|q\|_{L^2(\gamma)}
\|\operatorname{Tr}_\gamma\varphi\|_{L^2(\gamma)}
\le C_p\|q\|_{L^2(\gamma)}
\|\varphi\|_{W^{1,p'}(\Omega)}.
\]
Thus $q\delta_\gamma\in W^{-1,p}(\Omega)$.
Since $\Omega$ is a convex polyhedral domain, the Dirichlet
Laplacian defines an isomorphism
\[
-\Delta:W_0^{1,p}(\Omega)\longrightarrow W^{-1,p}(\Omega),
\qquad 1<p<\infty;
\]
see \cite[Corollary~2.7]{DrelichmanI_DuranR_OjeaI_2020}.
Consequently, for every $q\in L^2(\gamma)$, there exists a unique
state $u=S q\in W_0^{1,p}(\Omega)$ satisfying
\[
\|S q\|_{W^{1,p}(\Omega)}
 \le C\|q\delta_\gamma\|_{W^{-1,p}(\Omega)}
 \le C\|q\|_{L^2(\gamma)}.
\]
Since $p>6/5$, the embedding $W^{1,p}_0(\Om) \hookrightarrow L^2(\Om)$ gives 
 \begin{equation}\label{eq: cont of S from L2 to L2}
\|Sq\|_{L^2(\Om)}\le C\|q\|_{L^2(\gamma)}.
\end{equation}
 As a result, the operator $S$ is linear and continuous and we can introduce  the reduced functional $j \colon Q= L^2(\gamma) \to \R$ by 
\[
j(q):=J(q,u(q))=\frac{1}{2}\|Sq-u_d\|^2_{L^2(\Omega)}+\frac{\alpha}{2}\|q\|^2_{L^2(\gamma)}.
\]
Thus, we can rewrite the minimization problem \eqref{eq: minimize} as 
\begin{equation}\label{eq: minimize 2}
Minimize\quad j(q),\quad q\in Q_{ad}.
\end{equation}
\begin{theorem}[Existence and uniqueness of an optimal control]\label{thm:existence reduced}
There exists a unique solution $\bar{q}\in Q_{ad}$ to the minimization problem \eqref{eq: minimize 2}.
\end{theorem}
\begin{proof}
The proof is standard; see \cite[Thm.~10.1.6]{vexler2025numerical} for details. 
\end{proof}

The derivatives of $j(q)$ are  
\begin{subequations}
\begin{equation}\label{eq: j prime}
j'(q)(\delta q)=(u-u_d,\delta u)+\alpha\langle q,\delta q \rangle_\gamma
\end{equation}
and
\begin{equation}\label{eq: j double prime}
j''(q)(\delta q,\delta q)=\|\delta u\|^2_{L^2(\Om)}+\alpha\|\delta q\|^2_{L^2(\gamma)}
\end{equation}
\end{subequations}
where $u=Sq$ and $\delta u=S\delta q$.
We can represent $j'$ in terms of an adjoint variable $z=z(q)=z(u(q))\in H^2(\Om)\cap H^1_0(\Om)$, which is the solution to the equation 
\begin{equation}\label{eq: cont adjoint}
(\na z,\na \varphi)=({u}-{u}_d,\varphi),\quad \forall \varphi \in H^1_0(\Om).
\end{equation}
Thus, \eqref{eq: j prime} can be rewritten as
\begin{equation}\label{eq: j prime2}
\begin{aligned}
j'(q)(\delta q)&=(\na z,\na \delta u)+\alpha\langle q,\delta q \rangle_\gamma\\
&=\langle z, \delta q\rangle_\gamma+\alpha\langle q,\delta q \rangle_\gamma=\langle z+\alpha q,\delta q \rangle_\gamma.
\end{aligned}
\end{equation}
The  expression $(\na z,\na \delta u)$ requires a justification because  $\delta u$ generally does not belong to $H^1_0(\Om)$, but only to $W^{1,p}_0(\Om)$ for $6/5<p<2$. 
Because $\Omega$ is convex and $u-u_d\in L^2(\Omega)$, the adjoint
satisfies
$z\in H^2(\Omega)\cap H_0^1(\Omega).$
Consequently,
\[
z\in W^{1,6}(\Omega)\cap H_0^1(\Omega)
   \hookrightarrow W_0^{1,p'}(\Omega),
\qquad 2<p'<6.
\]
Moreover, $\delta u\in W_0^{1,p}(\Omega)\hookrightarrow L^2(\Omega)$.
Approximating $\delta u$ in $W_0^{1,p}(\Omega)$ by functions in
$C_0^\infty(\Omega)$ and using $-\Delta z=u-u_d$, we obtain
\[
(u-u_d,\delta u)_\Omega
   =(\nabla z,\nabla\delta u)_\Omega.
\]
Since $z\in W_0^{1,p'}(\Omega)$ is an admissible test function in the
linearized state equation,
\[
(\nabla\delta u,\nabla z)_\Omega
   =\langle\delta q,z\rangle_\gamma.
\]
Therefore,
\[
j'(q)(\delta q)
 =\langle z+\alpha q,\delta q\rangle_\gamma.
\]
The following characterization is standard for differentiable
convex minimization problems over convex admissible sets;
see, e.g., \cite[Lemma~5.2.1]{vexler2025numerical}.

\begin{lemma}\label{lemma: characterization of q}
A control $\bar{q}\in Q_{ad}$ is the solution to the problem \eqref{eq: minimize 2} if and only if
\[
j'(\bar{q})(\delta q-\bar{q})\geq 0,\quad \forall \delta q\in Q_{ad}.
\]
\end{lemma}

\subsection{Optimality conditions}

Using Lemma \ref{lemma: characterization of q}, the optimal solution $(\bar{q},\bar{u},\bar{z})$ satisfies:
\begin{subequations}

the state equation: $\bar{u}\in W^{1,p}_0(\Om)$  for $6/5<p<2$,
\begin{equation}\label{eq: state2}
(\na \bar{u},\na \varphi) = \langle \bar{q},\varphi\rangle_\gamma\quad \forall\varphi \in W^{1,p'}_0(\Om),\quad \frac{1}{p}+\frac{1}{p'}=1,
\end{equation}

the adjoint equation: $\bar{z}\in H^{1}_0(\Om)\cap H^2(\Om)$ 
\begin{equation}\label{eq: adjoint2}
(\na \bar{z},\na \varphi)=(\bar{u}-{u}_d,\varphi),\quad \forall \varphi \in H^1_0(\Om).
\end{equation}
and 
the optimality condition: $\bar{q}\in Q_{ad}$ 
\begin{equation}\label{eq: grad2}
 \langle\alpha \bar{q} +\bar{z}\mid_\gamma,\delta q-\bar{q}\rangle_\gamma\geq 0,\quad \forall \delta q\in Q_{ad}.
\end{equation}
\end{subequations}

\subsection{Regularity of the optimal solution}
We next establish regularity estimates for the optimal
solution. For this purpose, we use the Green's function representation.

Let $G(x,\xi)$ denote the Dirichlet Green's function
for the Laplacian on $\Omega$.
For every suitable $f$, the solution of
$-\Delta u=f$ in $\Omega$ with homogeneous Dirichlet
boundary conditions admits the representation
\begin{equation}\label{eq: Greens function representation}
	u(x)=\int_\Omega G(x,\xi)f(\xi)\,d\xi
	\qquad\text{for a.e. }x\in\Omega.
\end{equation}
The corresponding Green-potential representation for
the curve-supported source gives
\[
\bar u(x)
=\int_\gamma G(x,\xi)\bar q(\xi)\,d\mathcal H^1(\xi)
\qquad\text{for a.e. }x\in\Omega\setminus\gamma.
\]

Comparison with the fundamental solution in $\mathbb R^3$
by the maximum principle yields
\begin{equation}\label{Grees estimates}
	0\le G(x,\xi)\le \frac{1}{4\pi|x-\xi|}
	\qquad\text{for }x,\xi\in\Omega,\quad x\ne\xi.
\end{equation}
Moreover, convexity of $\Omega$ implies the gradient estimate
(cf. \cite{Fromm1993Potential})
\begin{equation}\label{eq: green's}
	|\nabla_x G(x,\xi)|
	\le \frac{C}{|x-\xi|^2}
	\qquad\text{for }x,\xi\in\Omega,\quad x\ne\xi.
\end{equation}

The following theorem provides regularity estimates for
the optimal control $\bar q$, state $\bar u$, and adjoint
state $\bar z$.
The estimates for the state and adjoint state are analogous
to those for pointwise controls in two dimensions;
see \cite[Theorem~10.2.7]{vexler2025numerical}.
In both settings, the source is supported on a set
of codimension two.

\begin{theorem}\label{thm: regularity}
Let $\Omega$ be a bounded convex polyhedral domain and
let $u_d\in L^\infty(\Omega)$. There exists $p_0=p_0(\Om)>2$ such that the optimal state $\bar u$, adjoint $\bar z$ and optimal control $\bar q$ satisfy the following global regularity properties:
\begin{subequations}
\begin{equation}\label{eq: state reg}
\bar{u}\in W^{1,p}_0(\Om)\cap L^r(\Om)\quad \forall p\in(1,2)\ \text{and}\ \forall r\in [1,\infty),
\end{equation}
\begin{equation}\label{eq: adjoint reg}
\bar{z}\in W^{2,p_0}(\Om)\cap W^{1,\infty}(\Om)\cap H^1_0(\Om),
\end{equation}
\begin{equation}\label{eq: control reg}
\bar{q}\in W^{1,\infty}(\gamma) .
\end{equation}
Moreover, for every open subdomain
$\Omega_0\Subset\Omega$ containing $\gamma$, we have
\begin{equation}\label{eq: adjoint reg2}
\bar{z}\in W^{2,r}(\Om_0),\quad \forall r\in [1,\infty)
\end{equation}
\end{subequations}
with the corresponding a priori estimates
\begin{subequations}
\begin{equation}\label{eq: state apriori2}
\|\bar{u}\|_{L^{r}(\Om)}\le Cr\|\bar{q}\|_{L^\infty(\gamma)},\quad \forall r\in [1,\infty),
\end{equation}
\begin{equation}\label{eq: state apriori}
\|\bar{u}\|_{W^{1,p}(\Om)}\le \frac{C}{(2-p)^\half}\|\bar{q}\|_{L^\infty(\gamma)},\quad 1<p<2
\end{equation}
\begin{equation}\label{eq: adjoint apriori}
\|\bar{z}\|_{H^{2}(\Om)}\le {C}\left(\|\bar{q}\|_{L^\infty(\gamma)}+\|u_d\|_{L^2(\Om)}\right),
\end{equation}
\begin{equation}\label{eq: adjoint apriori2}
	\|\bar z\|_{W^{2,r}(\Omega_0)}
	\le Cr^2\left(
	\|\bar q\|_{L^\infty(\gamma)}
	+\|u_d\|_{L^\infty(\Omega)}
	\right),
	\qquad 1\le r<\infty.
\end{equation}
\end{subequations}
The constants $C$ are independent of $p$ and $r$.
They may depend on $\Omega$ and the geometry of $\gamma$,
and, for the local estimate, on $\Omega_0$ and its
distance from $\partial\Omega$.
\end{theorem}
\begin{proof}
First, note that the box constraints imply that $\bar{q}\in L^\infty(\gamma)$ with $\|\bar{q}\|_{L^\infty(\gamma)}\le\max{\{|q_a|,|q_b|\}}$.
 
The state equation and boundedness of $\Omega$ imply
$\bar u\in W^{1,p}_0(\Omega)$ for every $1<p<2$.
By Sobolev embedding, $\bar u\in L^4(\Omega)$.
Since $u_d\in L^\infty(\Omega)$, the right-hand side
$\bar u-u_d$ of the adjoint equation belongs to $L^4(\Omega)$.
Hence, \cite[Theorem~2.2.6 and Corollary~2.2.16]
{vexler2025numerical} yield
\[
\bar z\in W^{2,p_0}(\Omega)\cap W^{1,\infty}(\Omega)
\]
for some $p_0=p_0(\Omega)>2$.
Since $\bar{z}\in W^{1,\infty}(\Om)$, it possesses a Lipschitz representative. Therefore $\bar{z}\circ \Gamma$
is Lipschitz on $[0,L]$, and hence 
  $\bar{z}|_\gamma\in W^{1,\infty}(\gamma)$. Using $\bar{q}=P_{[q_a,q_b]}\left(-\frac{1}{\alpha}\bar{z}\mid_\gamma\right)$, we also have   $\bar{q}\in W^{1,\infty}(\gamma)$. 

To establish the a priori estimates \eqref{eq: state apriori2} and \eqref{eq: state apriori}, we use again the Green's function representation of the solution 
\[
\bar{u}({x})
=\int_\gamma G({x},y)\bar{q}(y)\ d \Hg(y)
\]
and thus
\[
\abs{\bar u(x)} \le \|\bar{q}\|_{L^\infty(\gamma)} \int_\gamma G({x},y)\ d \Hg(y).
\]
Using the Green's function estimate \eqref{Grees estimates} and estimate \eqref{eq1:curve-kernel-one} from \cref{lem:curve-kernel-estimates}, we obtain
\begin{equation}\label{eq:state-pointwise}
|\bar{u}({x})|\le C \|\bar{q}\|_{L^\infty(\gamma)} \int_\gamma \frac{1}{\abs{x-y}}\ d \Hg(y) \le  C  \bigl(1+\abs{\ln d_\gamma(x)}\bigr) \|\bar{q}\|_{L^\infty(\gamma)},
\end{equation}
where $d_\gamma(x) = \operatorname{dist}(x,\gamma)$. Similarly, we get using the estimate for the gradient of the Green's function \eqref{eq: green's} and \eqref{eq1:curve-kernel-two} from \cref{lem:curve-kernel-estimates}
\begin{equation}\label{eq:state-gradient-pointwise}
|\nabla \bar{u}({x})|\le C \|\bar{q}\|_{L^\infty(\gamma)} \int_\gamma \frac{1}{\abs{x-y}^2}\ d \Hg(y) \le  \frac{C}{d_\gamma(x)} \|\bar{q}\|_{L^\infty(\gamma)}.
\end{equation}

We fix
\[
0<\varepsilon<
\min\{1,r_0,\operatorname{dist}(\gamma,\partial\Omega)\},
\]
where $r_0$ is the radius from \cref{lem:curve-tubular-integration}, and consider the $\varepsilon$-tubular neighborhood $T_\varepsilon(\gamma)$ of $\gamma$, see \eqref{eq:dist_and_tab_neigbourhood}. On the complement of this $\varepsilon$-tubular set, namely when
$d_\gamma(x)\ge\varepsilon$, the same Green's function estimates give
\[
 |\bar u(x)|+|\nabla\bar u(x)|
 \le C_\varepsilon\|\bar q\|_{L^\infty(\gamma)}.
\]
Thus, only the integrability of the right-hand sides of
\eqref{eq:state-pointwise} and
\eqref{eq:state-gradient-pointwise} inside
$T_\varepsilon(\gamma)$ remains to be considered.

Using \cref{lem:curve-tubular-integration} and
\eqref{eq:state-pointwise}, we obtain
\[
\|\bar u\|_{L^r(T_\varepsilon(\gamma))}^r
\le
C^r\|\bar q\|_{L^\infty(\gamma)}^r
\int_0^\varepsilon
(1+|\ln\rho|)^r\rho\,d\rho.
\]
As in the proof of \cite[Theorem~2.3.5]{vexler2025numerical},
we estimate this integral using
\[
\max_{0<\rho\le1}\rho(1-\ln\rho)^r
= r^r e^{1-r},
\qquad r\ge1.
\]
Indeed, the maximum is attained at $\rho=e^{1-r}$.
Since $\varepsilon\le1$, this yields
\[
\left(
\int_0^\varepsilon (1+|\ln\rho|)^r\rho\,d\rho
\right)^{1/r}
\le
\left(\varepsilon r^r e^{1-r}\right)^{1/r}
\le r.
\]
Combining this bound with the estimate on
$\Omega\setminus T_\varepsilon(\gamma)$ proves
\eqref{eq: state apriori2}, with a constant independent of $r$.

Similarly, \cref{lem:curve-tubular-integration} and
\eqref{eq:state-gradient-pointwise} imply
\[
\|\nabla\bar u\|_{L^p(T_\varepsilon(\gamma))}^p
\le
C^p\|\bar q\|_{L^\infty(\gamma)}^p
\int_0^\varepsilon \rho^{1-p}\,d\rho
=
C^p\|\bar q\|_{L^\infty(\gamma)}^p
\frac{\varepsilon^{2-p}}{2-p}.
\]
Since $\varepsilon\le1$, taking the $p$-th root and
including the estimate on
$\Omega\setminus T_\varepsilon(\gamma)$ yields
\[
\|\nabla\bar u\|_{L^p(\Omega)}
\le C(2-p)^{-1/p}\|\bar q\|_{L^\infty(\gamma)}.
\]
Moreover,
\[
(2-p)^{-1/p}
=
(2-p)^{-1/2}
(2-p)^{-\frac{2-p}{2p}}
\le C(2-p)^{-1/2},
\qquad 1<p<2,
\]
because the second factor is uniformly bounded.
Together with \eqref{eq: state apriori2}, applied with $r=p$,
this proves
\[
\|\bar u\|_{W^{1,p}(\Omega)}
\le \frac{C}{\sqrt{2-p}}\,
\|\bar q\|_{L^\infty(\gamma)},
\qquad 1<p<2,
\]
with a constant independent of $p$.

The inequality \eqref{eq: adjoint apriori} follows immediately from $H^2$ regularity, the triangle inequality and \eqref{eq: state apriori2},  
\[
\|\bar{z}\|_{H^2(\Om)}\le C\|\bar{u}-u_d\|_{L^2(\Om)}\le  C\left(\|\bar{u}\|_{L^2(\Om)}+\|u_d\|_{L^2(\Om)}\right)\le C\left(\|\bar{q}\|_{L^\infty(\gamma)}+\|u_d\|_{L^2(\Om)}\right).
\]
It remains to show  \eqref{eq: adjoint apriori2}.
Let $\Om_0\Subset \Om_1\Subset \Om$. 
For $r>2$, we apply the elliptic interior a priori regularity result (cf. \cite[Thm.~2.2.18]{vexler2025numerical}),  
\[
\|\bar{z}\|_{W^{2,r}(\Om_0)}\le Cr\left(\|\bar{u}-u_d\|_{L^r(\Om_1)}+\|\bar{u}-u_d\|_{L^2(\Om)} \right)\le Cr\|\bar{u}-u_d\|_{L^r(\Om)},
\]
where the constant $C$ is independent of $r$.
Using the triangle inequality and \eqref{eq: state apriori2},
\[
\|\bar{u}-u_d\|_{L^r(\Om)}\le \|\bar{u}\|_{L^r(\Om)}+\|u_d\|_{L^r(\Om)}\le Cr\|\bar{q}\|_{L^\infty(\gamma)}+C\|u_d\|_{L^\infty(\Om)}.
\]
Thus, we  obtain
\[
\|\bar{z}\|_{W^{2,r}(\Om_0)} \le Cr^2\|\bar{q}\|_{L^\infty(\gamma)}+Cr\|u_d\|_{L^\infty(\Om)}\le Cr^2\left(\|\bar{q}\|_{L^\infty(\gamma)}+\|u_d\|_{L^\infty(\Om)}\right).
\]
The case $1\le r\le 2$ follows by $H^2$-regularity and embedding $H^2(\Omega)\hookrightarrow W^{2,r}(\Om)$. 
\end{proof}

\section{Finite element approximation}\label{sec: fe discretization}

Let $\{\mathcal{T}_h\}_h$  be  a quasi-uniform shape-regular 
 family of conforming tetrahedral triangulations of $\Om$. For $\tau\in \mathcal{T}_h$, let $h_\tau=\operatorname{diam}(\tau)$ and set
 \[
 h=\max_{\tau} h_\tau.
 \]
 Then
\[
h_\tau\le h \le C |\tau|^{\frac{1}{3}}, \quad \forall \tau\in \mathcal{T}_h.
\]
Let 
\[
V_h=\{v_h\in H^1_0(\Omega)\ :\ v_h\mid_\tau\in \mathbb{P}_1\ \forall\tau\in \mathcal{T}_h\},
\] i.e. $V_h$ is the usual space of piecewise linear Lagrange finite elements.
For the finite element space $V_h$, we have the standard inverse inequality
\begin{equation}\label{eq: inverse}
\|\chi\|_{W^{l,r}(\Om)}\le Ch^{-l+3(\frac{1}{r}-\frac{1}{s})}\|\chi\|_{L^s(\Om)},\quad l=0,1,\quad 1\le s\le r\le \infty. 
\end{equation}
We also require the Ritz projection $R_h \colon W^{1,p}_0(\Omega) \to V_h$, $1\le p\le \infty$, defined by
\begin{equation}\label{eq:Ritz_proj}
(\nabla R_hv,\nabla \chi) = (\nabla v,\nabla \chi), \quad \forall \chi\in V_h.
\end{equation}
On convex polyhedral domains, the Ritz projection is stable in $W^{1,p}$ norm for any $p\in (1,\infty]$ (cf. \cite{DemlowA_LeykekhmanD_SchatzAH_WahlbinLB_2010a}), i.e.
\begin{equation}\label{eq: Stability Ritz}
\|\na R_hv\|_{L^p(\Om)}\le C\|\na v\|_{L^p(\Om)},\quad  1<p\le \infty.
\end{equation}
In the analysis we will use a nodewise interpolant $I_h \colon C(\bar{\Om})\cap H^1_0(\Omega) \to V_h$ with the usual approximation properties (cf., e.\,g.,~\cite[Theorem~3.1.5]{PGCiarlet_1978a})
\begin{equation}\label{eq:I_h approximation}
\|v-I_hv\|_{W^{1,r}(\Om)}\le Ch^{1+3(\frac{1}{r}-\frac{1}{2})}\|v\|_{H^{2}(\Om)}, \quad \text{for \ } 2\le r \le 6.
\end{equation}

By \cref{eq:curve-length-upper}, for all sufficiently small $h$, we immediately get 
\begin{equation}\label{eq:element-length}
\mathcal H^1(\gamma\cap\tau)\le Ch
\qquad\text{for every }\tau\in\mathcal T_h.
\end{equation}

\subsection{The weight function}
Recall that
\begin{equation}
    \label{eq: distatnce to curve}
    d_\gamma(x) :=\operatorname{dist}(x,\gamma)=\inf_{y\in \gamma}|x-y|
\end{equation}
and introduce the positive weight function $\sigma_h$,
\begin{equation}\label{eq: sigma weight}
\sigma_h(x) :=\sqrt{d_\gamma^2(x)+h^2}.
\end{equation}
The weight function $\sigma_h$ has the following properties:
\begin{subequations}
\begin{align}
0<h\le \sigma_h&\le \operatorname{diam}(\Omega)+ h, \label{eq: property 1 of sigma}\\
|\na \sigma_h(x)|&\le 1,\quad \text{a.e. $x\in \Om$,}\label{eq: property 2 of sigma}\\
\|\sigma_h^{-1}\|_{L^2(\Om)}&\le C\ell_h^{1/2}, \label{eq: property  of sigma in L2}\\
\max_\tau{\sigma_h}&\le C\min_\tau{\sigma_h}, \quad \forall \tau \in \,\mathcal{T}_h\label{eq: property 4 of sigma},
\end{align}
\end{subequations}
where in \eqref{eq: property  of sigma in L2}  and for the rest of the paper, we define 
\begin{equation}\label{eq: ell_h}
    \ell_h = (1+\abs{\ln{h}}).
\end{equation}
The above properties are not difficult to establish. Thus, \eqref{eq: property 1 of sigma} is obvious and the property \eqref{eq: property 2 of sigma} follows from the observation
\[
|\na \sigma_h|=\frac{d_\gamma}{\sigma_h}|\na d_\gamma|\le 1.
\]
The property \eqref{eq: property 4 of sigma} is a simple consequence of \eqref{eq: property 2 of sigma} and
\[
\max_{x\in\tau}{\sigma_h}(x)\le \min_{y\in\tau}{\sigma_h(y)}+\operatorname{diam}(\tau)\le \min_{y\in\tau}{\sigma_h(y)}+h\le 2\min_{y\in\tau}{\sigma_h(y)}.
\]

To prove \eqref{eq: property  of sigma in L2}, fix
$0<\varepsilon<\min\{r_0,\operatorname{dist}(\gamma,\partial\Omega)\}$,
where $r_0$ is the radius from
\cref{lem:curve-tubular-integration}.
Applying that lemma with $F(\rho)=(\rho^2+h^2)^{-1}$ gives
\[
\begin{aligned}
	\|\sigma_h^{-1}\|_{L^2(\Omega)}^2
	&=
	\int_{T_\varepsilon(\gamma)}
	\frac{dx}{d_\gamma(x)^2+h^2}
	+
	\int_{\Omega\setminus T_\varepsilon(\gamma)}
	\frac{dx}{d_\gamma(x)^2+h^2}\\
	&\le
	C\int_0^\varepsilon\frac{\rho}{\rho^2+h^2}\,d\rho
	+\varepsilon^{-2}|\Omega|\\
	&=
	\frac{C}{2}\ln\left(1+\frac{\varepsilon^2}{h^2}\right)
	+\varepsilon^{-2}|\Omega|
	\le C\ell_h.
\end{aligned}
\]
Taking square roots proves
\eqref{eq: property  of sigma in L2}.

\subsection{Weighted results}
Similarly to \cite[p.~998]{ApelT_BenedixO_SirchD_VexlerB_2011a} and \cite[Lem.~2.5]{LeykekhmanD_VexlerB_2016b} we will show a weighted Hardy inequality. We start with a one-dimensional radial inequality.
\begin{lemma}[Radial weighted Hardy inequality]
	\label{lem:radial-weighted-hardy}
	Let $R>0$ and $0<h\le R$. Let
	$v\in H^1((0,R))$ satisfy $v(R)=0$.
	Then
	\[
	\int_0^R
	\frac{|v(\rho)|^2}{\rho^2+h^2}\,\rho\,d\rho
	\le
	C \left(1+\ln\frac{R}{h}\right)^2
	\int_0^R |v'(\rho)|^2\rho\,d\rho.
	\]
\end{lemma}
\begin{proof}
	Set $\sigma(\rho):=\sqrt{\rho^2+h^2}$ and
	$L_h:=1+\ln(R/h)$. Then
	\[
	\frac{d}{d\rho}
	\bigl(\ln\sigma(\rho)-\ln\sigma(0)\bigr)
	=\frac{\rho}{\sigma(\rho)^2}.
	\]
	Moreover,
	\[
	\frac{\ln\sigma(\rho)-\ln\sigma(0)}{\rho}
	\le \frac{C L_h}{\sigma(\rho)},
	\qquad 0<\rho\le R.
	\]
	Indeed, for $0<\rho\le h$, we have
	\[
	\ln\sigma(\rho)-\ln\sigma(0)
	=\frac12\ln\left(1+\frac{\rho^2}{h^2}\right)
	\le \frac{\rho^2}{2h^2},
	\]
	and hence
	\[
	\frac{\ln\sigma(\rho)-\ln\sigma(0)}{\rho}
	\le \frac{1}{2h}
	\le \frac{C}{\sigma(\rho)}.
	\]
	For $h\le\rho\le R$, we use
	$\sigma(\rho)\le\sqrt{2}\rho$ to obtain
	\[
	\frac{\ln\sigma(\rho)-\ln\sigma(0)}{\rho}
	= \frac{1}{2\rho} 
	 \ln\left(1+\frac{\rho^2}{h^2}\right)
	\le \frac{1}{2\rho}\ln\left(1+\frac{R^2}{h^2}\right)
	\le \frac{C L_h}{\sigma(\rho)}.
	\]
	
	Since $v(R)=0$ and
	$\ln\sigma(\rho)-\ln\sigma(0)$ vanishes at $\rho=0$,
	integration by parts gives
	\[
	\begin{aligned}
		\int_0^R
		\frac{|v(\rho)|^2}{\sigma(\rho)^2}\rho\,d\rho
		&=
		-2\int_0^R
		\bigl(\ln\sigma(\rho)-\ln\sigma(0)\bigr)
		v(\rho)v'(\rho)\,d\rho\\
		&\le
		C L_h\int_0^R
		\frac{|v(\rho)|}{\sigma(\rho)}
		|v'(\rho)|\rho\,d\rho\\
		&\le
		C L_h
		\left(
		\int_0^R\frac{|v(\rho)|^2}{\sigma(\rho)^2}
		\rho\,d\rho
		\right)^{1/2}
		\left(
		\int_0^R|v'(\rho)|^2\rho\,d\rho
		\right)^{1/2}.
	\end{aligned}
	\]
	Dividing by the first factor on the right, unless it
	vanishes, and squaring proves the assertion.
\end{proof}

\begin{lemma}[Weighted Hardy inequality]
	\label{lem:weighted-hardy}
	Let $\gamma\Subset\Omega$ be a compact embedded $C^{1,1}$
	curve with two endpoints, and let
	\[
	\sigma_h(x):=\sqrt{d_\gamma(x)^2+h^2},
	\qquad
	\ell_h:=1+|\ln h|.
	\]
	There exist constants $C>0$ and $h_0>0$, depending only
	on $\Omega$ and $\gamma$, such that
	\[
	\|\sigma_h^{-1}v\|_{L^2(\Omega)}
	\le C\ell_h\|\nabla v\|_{L^2(\Omega)}
	\]
	for every $v\in H^1_0(\Omega)$ and $0<h\le h_0$.
\end{lemma}
\begin{proof}
	By density, it suffices to consider $v\in C_c^\infty(\Omega)$.
	Let $\Phi$, $J$, $R$, $r$, and $W^0$ be as in the
	tubular coordinate construction in
	\cref{lem:curve-tubular-coordinates}. Choose
	$\eta\in C_c^\infty(W^0)$ with $0\le\eta\le1$ and
	$\eta=1$ on $T_r(\gamma)$, and set $w:=\eta v$.
	
	Since $(1-\eta)v$ vanishes on $T_r(\gamma)$,
	the Poincaré inequality gives
	\[
	\|\sigma_h^{-1}(1-\eta)v\|_{L^2(\Omega)}
	\le r^{-1}\|v\|_{L^2(\Omega)}
	\le C\|\nabla v\|_{L^2(\Omega)}.
	\]
	For the remaining term, set
	$V(s,\rho,\theta):=w(\Phi(s,\rho,\theta))$.

Since $\operatorname{supp}w\Subset W^0$, the change of
variables $x=\Phi(s,\rho,\theta)$ gives
\[
\begin{aligned}
	\|\sigma_h^{-1}w\|_{L^2(\Omega)}^2
	&=
	\int_J\int_0^{2\pi}\int_0^R
	\frac{|V(s,\rho,\theta)|^2}
	{d_\gamma(\Phi(s,\rho,\theta))^2+h^2}
	|\det D\Phi(s,\rho,\theta)|
	\,d\rho\,d\theta\,ds\\
	&\le
	C\int_J\int_0^{2\pi}\int_0^R
	\frac{|V(s,\rho,\theta)|^2}{\rho^2+h^2}
	\,\rho\,d\rho\,d\theta\,ds,
\end{aligned}
\]
where we used $d_\gamma(\Phi(s,\rho,\theta))\ge\rho$
and the upper Jacobian bound from \cref{lem:curve-tubular-coordinates}.

For fixed $s$ and $\theta$, the function
$V(s,\cdot,\theta)$ belongs to $H^1(0,R)$ and vanishes
near $\rho=R$. Thus, \cref{lem:radial-weighted-hardy}
and $1+\ln(R/h)\le C\ell_h$ imply
\[
\|\sigma_h^{-1}w\|_{L^2(\Omega)}^2
\le C\ell_h^2
\int_J\int_0^{2\pi}\int_0^R
|\partial_\rho V(s,\rho,\theta)|^2
\,\rho\,d\rho\,d\theta\,ds.
\]

Finally, the chain rule and $|\partial_\rho\Phi|=1$ give
\[
|\partial_\rho V(s,\rho,\theta)|
=
\bigl|\nabla w(\Phi(s,\rho,\theta))
\cdot\partial_\rho\Phi(s,\rho,\theta)\bigr|
\le |\nabla w(\Phi(s,\rho,\theta))|.
\]
Using the lower Jacobian bound from \cref{lem:curve-tubular-coordinates} and transforming back
to $W^0$, we obtain
\[
\begin{aligned}
	&\int_J\int_0^{2\pi}\int_0^R
	|\partial_\rho V(s,\rho,\theta)|^2
	\,\rho\,d\rho\,d\theta\,ds\\
	&\qquad\le
	C\int_J\int_0^{2\pi}\int_0^R
	|\nabla w(\Phi(s,\rho,\theta))|^2
	|\det D\Phi(s,\rho,\theta)|
	\,d\rho\,d\theta\,ds\\
	&\qquad=
	C\|\nabla w\|_{L^2(W^0)}^2
	=
	C\|\nabla w\|_{L^2(\Omega)}^2.
\end{aligned}
\]
	Since $\eta$ is fixed, the Poincaré inequality further yields
	\[
	\|\nabla w\|_{L^2(\Omega)}
	\le C\bigl(\|\nabla v\|_{L^2(\Omega)}
	+\|v\|_{L^2(\Omega)}\bigr)
	\le C\|\nabla v\|_{L^2(\Omega)}.
	\]
	Combining the two estimates proves the assertion for
	smooth $v$. For each fixed $h>0$, we have
	$\sigma_h^{-1}\le h^{-1}$, so passage to the limit
	extends the estimate to every $v\in H^1_0(\Omega)$.
\end{proof}

The following result is a discrete Sobolev-type inequality
for the restriction of finite element functions to $\gamma$.
It is analogous to the classical discrete Sobolev inequality
in two dimensions; see
\cite[Lemma~4.9.2]{BrennerScott2008}. The proof follows the radial argument of
\cite[Lemma~4.9.2]{BrennerScott2008}, combined with
the tubular coordinates from \cref{lem:curve-tubular-coordinates}
and \cref{cor:curve-ball-integral}.

\begin{lemma}[Discrete Sobolev-type inequality]
	\label{lem:discrete-sobolev}
	Let $\{\mathcal T_h\}_h$ be a quasi-uniform, shape-regular
	family of tetrahedral meshes, and let
	$V_h\subset H^1_0(\Omega)$ be the associated space of
	continuous piecewise linear functions.
	There exist constants $C>0$ and $h_0>0$ such that
	\[
	\|v_h\|_{L^2(\gamma)}
	\le C\ell_h^{\half}\|\nabla v_h\|_{L^2(\Omega)}
	\qquad
	\text{for all }v_h\in V_h,\quad 0<h\le h_0.
	\]
	The constants depend only on $\Omega$, the geometry of
	$\gamma$, and the mesh regularity constants.
\end{lemma}
\begin{proof}
	Let $\Phi$, $J$, $R$, $r$, and $W^0$ be as in
	\cref{lem:curve-tubular-coordinates}, with $\Gamma:[0,L]\to\gamma$
	the arclength parametrization of the original curve.
	Choose $\eta\in C_c^\infty(W^0)$ such that
	$\eta=1$ on $T_r(\gamma)$, and set $w:=\eta v_h$.
	We take $h_0>0$ sufficiently small that
	\[
	h_0<\min\{r,R,1\},
	\qquad
	2h_0\le r_{\rm len},
	\qquad
	2h_0<\operatorname{dist}(\gamma,\partial\Omega).
	\]
	
	First, for $x\in\gamma$, consider the patch
	\[
	\omega_{x,h}
	:=
	\bigcup\{\tau\in\mathcal T_h:\tau\cap B_h(x)\ne\varnothing\}.
	\]
	Since $\operatorname{diam}(\tau)\le h$, this patch is
	contained in $B_{2h}(x)$. The elementwise inverse
	inequality and quasi-uniformity therefore give
	\[
	\|\nabla v_h\|_{L^\infty(B_h(x))}
	\le Ch^{-3/2}
	\|\nabla v_h\|_{L^2(B_{2h}(x))}.
	\]
	
	For $s\in[0,L]$, set $x=\Gamma(s)$ and
	\[
	V(s,\rho,\theta):=w(\Phi(s,\rho,\theta)).
	\]
	Then $V(s,0,\theta)=v_h(x)$, and $V$ vanishes near
	$\rho=R$ because $w$ has compact support in $W^0$.
	The fundamental theorem of calculus yields
	\[
	v_h(x)
	=
	-\int_0^h\partial_\rho V(s,\rho,\theta)\,d\rho
	-\int_h^R\partial_\rho V(s,\rho,\theta)\,d\rho.
	\]
	
	For $0\le\rho\le h$, we have
	$|\Phi(s,\rho,\theta)-x|=\rho$ and $\eta=1$.
	Consequently, using $|\partial_\rho\Phi|=1$ and
	the local inverse estimate above,
	\[
	\begin{aligned}
		\left|
		\int_0^h\partial_\rho V(s,\rho,\theta)\,d\rho
		\right|^2
		&\le
		h^2\|\nabla v_h\|_{L^\infty(B_h(x))}^2\\
		&\le
		Ch^{-1}\|\nabla v_h\|_{L^2(B_{2h}(x))}^2.
	\end{aligned}
	\]
	For the outer integral, the Cauchy--Schwarz inequality gives
	\[
	\begin{aligned}
		\left|
		\int_h^R\partial_\rho V(s,\rho,\theta)\,d\rho
		\right|^2
		&\le
		\left(\int_h^R
		|\partial_\rho V(s,\rho,\theta)|^2\rho\,d\rho\right)
		\left(\int_h^R\frac{d\rho}{\rho}\right)\\
		&\le
		C\ell_h\int_h^R
		|\partial_\rho V(s,\rho,\theta)|^2\rho\,d\rho.
	\end{aligned}
	\]
	Combining these estimates and averaging over $\theta$,
	we obtain
	\[
	\begin{aligned}
		|v_h(\Gamma(s))|^2
		&\le
		Ch^{-1}
		\|\nabla v_h\|_{L^2(B_{2h}(\Gamma(s)))}^2\\
		&\quad+
		C\ell_h\int_0^{2\pi}\int_h^R
		|\partial_\rho V(s,\rho,\theta)|^2
		\rho\,d\rho\,d\theta.
	\end{aligned}
	\]
	
	We now integrate over $s\in[0,L]$.
	Since $\Gamma$ is parametrized by arclength,
	\cref{cor:curve-ball-integral} implies
	\[
	h^{-1}\int_0^L
	\|\nabla v_h\|_{L^2(B_{2h}(\Gamma(s)))}^2\,ds
	=
	h^{-1}\int_\gamma
	\|\nabla v_h\|_{L^2(B_{2h}(x))}^2
	\,d\mathcal H^1(x)
	\le C\|\nabla v_h\|_{L^2(\Omega)}^2.
	\]
	For the second term, we use
	$|\partial_\rho V|\le|\nabla w\circ\Phi|$
	and the lower Jacobian bound
	$|\det D\Phi|\ge c\rho$ to obtain
	\[
	\begin{aligned}
		&\int_0^L\int_0^{2\pi}\int_h^R
		|\partial_\rho V(s,\rho,\theta)|^2
		\rho\,d\rho\,d\theta\,ds\\
		&\qquad\le
		C\int_0^L\int_0^{2\pi}\int_h^R
		|\nabla w(\Phi(s,\rho,\theta))|^2
		|\det D\Phi(s,\rho,\theta)|
		\,d\rho\,d\theta\,ds\\
		&\qquad\le C\|\nabla w\|_{L^2(\Omega)}^2.
	\end{aligned}
	\]
	Finally, since $\eta$ is fixed, the Poincaré inequality gives
	\[
	\|\nabla w\|_{L^2(\Omega)}
	\le C\bigl(\|\nabla v_h\|_{L^2(\Omega)}
	+\|v_h\|_{L^2(\Omega)}\bigr)
	\le C\|\nabla v_h\|_{L^2(\Omega)}.
	\]
	It follows that
	\[
	\|v_h\|_{L^2(\gamma)}^2
	\le C\ell_h\|\nabla v_h\|_{L^2(\Omega)}^2,
	\]
	which proves the assertion.
\end{proof}

\subsection{Supplementary results for the state equation}
We consider the basic finite element result for the state equation
\begin{equation}\label{eq: state4}
\begin{aligned}
-\Delta {u}  &={q}\delta_\gamma \quad \text{in $\Omega$}\\
{u} &=0  \quad \text{on $\partial\Omega$}.
\end{aligned}
\end{equation}
The Ritz projection $u_h=R_hu$ satisfies
\begin{equation}\label{eq: fe state}
(\na u_h,\na \chi) =\langle q, \chi \rangle_\gamma\quad \forall
\chi\in V_h.  
\end{equation}
Naturally, we have the Galerkin orthogonality
\[
(\na (u-u_h),\na \chi) =0, \quad \forall
\chi\in V_h.
\]

We establish an $L^2$ error estimate for the finite element
approximation of the state equation with $q\in L^2(\gamma)$.

\begin{lemma}\label{lem:state-L2-error}
	Let $q\in L^2(\gamma)$, and let $u$ and $u_h\in V_h$
	be the solutions of \eqref{eq: state4} and
	\eqref{eq: fe state}, respectively.
	Then, for all sufficiently small $h$,
	\[
	\|u-u_h\|_{L^2(\Omega)}
	\le Ch\|q\|_{L^2(\gamma)},
	\]
	where $C$ is independent of $h$ and $q$, but may depend
	on $\Omega$, the geometry of $\gamma$, and the
	mesh regularity constants.
\end{lemma}
\begin{proof}
	Set $e:=u-u_h$, and let $z\in H^2(\Omega)\cap H^1_0(\Omega)$
	solve the dual problem
	\[
	-\Delta z=e \quad\text{in }\Omega,
	\qquad
	z=0 \quad\text{on }\partial\Omega.
	\]
	By convexity of $\Omega$,
	\[
	\|z\|_{H^2(\Omega)}\le C\|e\|_{L^2(\Omega)}.
	\]
	Let $R_hz\in V_h$ denote the Ritz projection of $z$,
	defined by
	\[
	(\nabla R_hz,\nabla\chi)
	=(\nabla z,\nabla\chi)
	\qquad\text{for every }\chi\in V_h.
	\]
	The continuous and discrete state equations imply
	\[
	(u,e)=\langle q,z\rangle_\gamma
	\]
	and
	\[
	(u_h,e)
	=(\nabla u_h,\nabla z)
	=(\nabla u_h,\nabla R_hz)
	=\langle q,R_hz\rangle_\gamma.
	\]
	Consequently,
	\[
	\|e\|_{L^2(\Omega)}^2
	=\langle q,z-R_hz\rangle_\gamma
	\le
	\|q\|_{L^2(\gamma)}
	\|z-R_hz\|_{L^2(\gamma)}.
	\]
	
	To estimate the last factor, let $I_hz\in V_h$ be the
	nodal interpolant of $z$ and write
	\[
	z-R_hz=(z-I_hz)+(I_hz-R_hz).
	\]
	Set
	\[
	\mathcal T_{\gamma,h}
	:=\{\tau\in\mathcal T_h:\tau\cap\gamma\ne\varnothing\}.
	\]
	The local length estimate \eqref{eq:curve-length-upper}
	implies $\mathcal H^1(\gamma\cap\tau)\le Ch$.
	Together with the local interpolation estimate in
	three dimensions, this gives
	\[
	\begin{aligned}
		\|z-I_hz\|_{L^2(\gamma)}^2
		&\le
		\sum_{\tau\in\mathcal T_{\gamma,h}}
		\mathcal H^1(\gamma\cap\tau)
		\|z-I_hz\|_{L^\infty(\tau)}^2\\
		&\le
		C\sum_{\tau\in\mathcal T_{\gamma,h}}
		h\,h_\tau |z|_{H^2(\tau)}^2
		\le Ch^2\|z\|_{H^2(\Omega)}^2.
	\end{aligned}
	\]
	Here we used
	$\|z-I_hz\|_{L^\infty(\tau)}
	\le Ch_\tau^{1/2}|z|_{H^2(\tau)}$.
	
	For any $\chi_h\in V_h$, the inverse inequality and
	quasi-uniformity similarly yield
	\[
	\begin{aligned}
		\|\chi_h\|_{L^2(\gamma)}^2
		&\le
		Ch\sum_{\tau\in\mathcal T_{\gamma,h}}
		\|\chi_h\|_{L^\infty(\tau)}^2\\
		&\le
		Ch^{-2}\sum_{\tau\in\mathcal T_{\gamma,h}}
		\|\chi_h\|_{L^2(\tau)}^2
		\le Ch^{-2}\|\chi_h\|_{L^2(\Omega)}^2.
	\end{aligned}
	\]
	Applying this estimate to $\chi_h=I_hz-R_hz$ and using
	the standard interpolation and Ritz projection estimates,
	we obtain
	\[
	\begin{aligned}
		\|I_hz-R_hz\|_{L^2(\gamma)}
		&\le Ch^{-1}\|I_hz-R_hz\|_{L^2(\Omega)}\\
		&\le Ch^{-1}\bigl(
		\|I_hz-z\|_{L^2(\Omega)}
		+\|z-R_hz\|_{L^2(\Omega)}\bigr)\\
		&\le Ch\|z\|_{H^2(\Omega)}.
	\end{aligned}
	\]
	Thus,
	\[
	\|z-R_hz\|_{L^2(\gamma)}
	\le Ch\|z\|_{H^2(\Omega)}.
	\]
	Inserting this into the duality identity gives
	\[
	\|e\|_{L^2(\Omega)}^2
	\le Ch\|q\|_{L^2(\gamma)}\|z\|_{H^2(\Omega)}
	\le Ch\|q\|_{L^2(\gamma)}\|e\|_{L^2(\Omega)},
	\]
	which proves the assertion.
\end{proof}

\subsection{Interior error estimate}
The following interior pointwise error estimate follows
from \cite[Corollary~5.1]{AHSchatz_LBWahlbin_1977a}.

\begin{theorem}\label{thm: Schatz_Wahlbin_1977}
	Let $1\le p\le\infty$.
	There exist constants $C>0$, $\kappa>0$, and $h_0>0$
	such that the following holds.
	Let $0<h\le h_0$, $x_0\in\Omega$, and $d\ge\kappa h$
	satisfy $B_{2d}(x_0)\Subset\Omega$.
	For $w\in H^1_0(\Omega)$ with
	$w\in W^{2,\infty}(B_d(x_0))$, let $w_h=R_hw$.
	Then
	\[
	|(w-w_h)(x_0)|
	\le
	Ch^2\ell_h|w|_{W^{2,\infty}(B_d(x_0))}
	+
	Cd^{-3/p}\|w-w_h\|_{L^p(B_d(x_0))}.
	\]
	The constants are independent of $h$, $x_0$, $d$, and $w$.
	For $p=\infty$, the factor $d^{-3/p}$ is understood as $1$.
\end{theorem}

\subsection{Weighted error estimate for the state equation}

\begin{theorem}\label{thm: weighted error in L2}
Let $u$ be the solution to \eqref{eq: state4} with $q\in L^2(\gamma)$ and $u_h=R_hu$ be its Galerkin projection. Then there exists a constant $C$ independent of $h$ such that
\[
\|\sigma_h(u-u_h)\|_{L^2(\Om)}\le Ch^2\ell_h^{3/2}\|q\|_{L^2(\gamma)}.
\]
\end{theorem}

\begin{proof}
	To estimate $\|\sigma_h(u-u_h)\|_{L^2(\Omega)}$, we use
	a dyadic tubular decomposition
	\[
	\Omega=\Omega_*\cup\bigcup_{j=j_0}^J\Omega_j,
	\]
	where $j_0\in\mathbb Z$ is the largest integer such that
	$d_{j_0}:=2^{-j_0}\ge\operatorname{diam}(\Omega)$, and
	\[
	\Omega_*=\{x\in\Omega:d_\gamma(x)\le d_{J+1}\},
	\qquad
	\Omega_j=\{x\in\Omega:d_{j+1}\le d_\gamma(x)\le d_j\}.
	\]
	Here $d_j:=2^{-j}$, $C_*$ is a sufficiently large
	constant, and $J$ is chosen such that
	\[
	\frac{C_*}{2}h\le d_J\le C_*h.
	\]
	Let
	\[
	d_{\gamma,\partial\Omega}:=\operatorname{dist}(\gamma,\partial\Omega),
	\qquad
	r_0:=\frac1{16}\min\{d_{\gamma,\partial\Omega},r_{\rm len}\},
	\]
	and, for every shell $\Omega_j$, define
	\[
	r_j:=\min\left\{\frac{d_j}{4},r_0\right\}.
	\]
	Choose $C_*$ sufficiently large and $h$ sufficiently
	small so that $r_j\ge\kappa h$ for every $j_0\le j\le J$,
	where $\kappa$ is the constant from
	\cref{thm: Schatz_Wahlbin_1977}.
	For every $x\in\gamma$, we then have
	\[
	B_{2r_j}(x)\Subset\Omega,
	\qquad
	B_{2r_j}(x)\cap\Omega_j=\varnothing.
	\]
	
	Set $e:=u-u_h$. Since $\sigma_h^2\simeq d_j^2$ on
	$\Omega_j$, we obtain
	\[
	\|\sigma_he\|_{L^2(\Omega)}^2
	\le
	C h^2\|e\|_{L^2(\Omega_*)}^2
	+
	C\sum_{j=j_0}^J d_j^2\|e\|_{L^2(\Omega_j)}^2.
	\]
	By \cref{lem:state-L2-error},
	\[
	\|e\|_{L^2(\Omega_*)}
	\le \|e\|_{L^2(\Omega)}
	\le Ch\|q\|_{L^2(\gamma)}.
	\]
	
	To estimate $\|e\|_{L^2(\Omega_j)}$, consider the dual
	problem
	\[
	-\Delta v=e\chi_{\Omega_j}\quad\text{in }\Omega,
	\qquad
	v=0\quad\text{on }\partial\Omega,
	\]
	and let $v_h=R_hv$.
	By convexity of $\Omega$,
	\[
	\|v\|_{H^2(\Omega)}
	\le C\|e\|_{L^2(\Omega_j)}.
	\]
	Moreover, $v$ is harmonic in $B_{2r_j}(x)$ for every
	$x\in\gamma$. The state equations and Galerkin
	orthogonality give
	\[
	\begin{aligned}
		\|e\|_{L^2(\Omega_j)}^2
		&=(e,e\chi_{\Omega_j})
		=(\nabla e,\nabla v)\\
		&=(\nabla u,\nabla(v-v_h))
		=\int_\gamma q(v-v_h)\,d\mathcal H^1\\
		&\le
		\|q\|_{L^2(\gamma)}\|v-v_h\|_{L^2(\gamma)}.
	\end{aligned}
	\]
	The terms involving $\nabla u$ are understood as
	$L^p$--$L^{p'}$ pairings, with any fixed $6/5<p<2$.
	
	Applying \cref{thm: Schatz_Wahlbin_1977} with $p=2$
	at every $x\in\gamma$, we obtain
	\[
	|(v-v_h)(x)|
	\le
	Ch^2\ell_h|v|_{W^{2,\infty}(B_{r_j}(x))}
	+
	Cr_j^{-3/2}\|v-v_h\|_{L^2(B_{r_j}(x))}.
	\]
	Squaring and integrating over $\gamma$ yields
	\[
	\begin{aligned}
		\|v-v_h\|_{L^2(\gamma)}^2
		&\le
		Ch^4\ell_h^2
		\int_\gamma
		|v|_{W^{2,\infty}(B_{r_j}(x))}^2
		\,d\mathcal H^1(x)\\
		&\quad+
		Cr_j^{-3}
		\int_\gamma
		\|v-v_h\|_{L^2(B_{r_j}(x))}^2
		\,d\mathcal H^1(x).
	\end{aligned}
	\]	
	For the second term, \cref{cor:curve-ball-integral} and the
	standard $L^2$ error estimate for the Ritz projection give
	\[
		\int_\gamma
		\|v-v_h\|_{L^2(B_{r_j}(x))}^2\,d\mathcal H^1(x)
		\le Cr_j\|v-v_h\|_{L^2(\Omega)}^2
		\le Cr_jh^4\|v\|_{H^2(\Omega)}^2.
	\]
	For the first term, harmonicity of $v$ in $B_{2r_j}(x)$
	and the interior estimate for harmonic functions imply
	\[
	|v|_{W^{2,\infty}(B_{r_j}(x))}
	\le
	Cr_j^{-3/2}\|D^2v\|_{L^2(B_{2r_j}(x))}.
	\]
	Applying \cref{cor:curve-ball-integral} once more, we obtain
	\[
		\int_\gamma
		|v|_{W^{2,\infty}(B_{r_j}(x))}^2\,d\mathcal H^1(x)
		\le
		Cr_j^{-3}\int_\gamma
		\|D^2v\|_{L^2(B_{2r_j}(x))}^2\,d\mathcal H^1(x)
		\le Cr_j^{-2}\|v\|_{H^2(\Omega)}^2.
	\]
	Combining these estimates yields
	\[
	\|v-v_h\|_{L^2(\gamma)}
	\le Ch^2\ell_h r_j^{-1}\|v\|_{H^2(\Omega)}
	\le Ch^2\ell_h r_j^{-1}\|e\|_{L^2(\Omega_j)}.
	\]
	Thus,
	\[
	\|e\|_{L^2(\Omega_j)}^2
	\le
	Ch^2\ell_h r_j^{-1}
	\|q\|_{L^2(\gamma)}\|e\|_{L^2(\Omega_j)},
	\]
	and consequently,
	\[
	\|e\|_{L^2(\Omega_j)}^2
	\le Ch^4\ell_h^2r_j^{-2}\|q\|_{L^2(\gamma)}^2.
	\]
	
	We have therefore shown that
	\[
	\|\sigma_he\|_{L^2(\Omega)}^2
	\le
	Ch^4\|q\|_{L^2(\gamma)}^2
	\left(
	1+\ell_h^2\sum_{j=j_0}^J
	\left(\frac{d_j}{r_j}\right)^2
	\right).
	\]
	If $d_j\le4r_0$, then $r_j=d_j/4$ and $d_j/r_j=4$.
	If $d_j>4r_0$, then $r_j=r_0$; the number of such
	indices and the ratios $d_j/r_0$ are bounded
	independently of $h$. Since $J-j_0+1\le C\ell_h$,
	we conclude that
	\[
	\sum_{j=j_0}^J\left(\frac{d_j}{r_j}\right)^2
	\le C\ell_h.
	\]
	Hence
	\[
	\|\sigma_he\|_{L^2(\Omega)}^2
	\le Ch^4\ell_h^3\|q\|_{L^2(\gamma)}^2.
	\]
	Taking square roots proves the assertion.
\end{proof}

The above result is of independent interest.
Using Theorem \ref{thm: weighted error in L2} and the definition of $\sigma_h$ \eqref{eq: sigma weight} and the property \eqref{eq: property  of sigma in L2}, 
we immediately obtain the following  results.
\begin{corollary}\label{cor: L1}
    Under the assumptions of Theorem \ref{thm: weighted error in L2},
    \[
\|u-u_h\|_{L^1(\Omega)}\le Ch^2\ell_h^{2}\|q\|_{L^2(\gamma)}
\]
and
    \[
\|u-u_h\|_{L^2(\Omega_d)}\le Ch^2\ell_h^{3/2}d^{-1}\|q\|_{L^2(\gamma)},
\]
with a constant $C$ independent of $h$ and $d$, where
\[
\Omega_d := \{ x\in \Omega\ :\ \operatorname{dist}(x,\gamma)> d\}.
\]
\end{corollary}

\section{Variational control discretization}\label{sec: variational}

In this section we analyze and derive improved error estimates in the case of variational discretization, meaning that the control $q$ is not discretized. Approximation of \(\gamma\) by a polygonal
curve would introduce an additional geometric variational crime and
is not considered here.

First we introduce a discrete solution operator $S_h:Q\to V_h$, which maps a given control $q\in Q$ to the finite element solution $u_h=u_h(q)=S_hq\in V_h$ by
\begin{equation}\label{eq: discrete state}
(\na u_h,\na \chi) = \int_\gamma q\chi \ d\Hg\quad \forall\chi \in V_h. 
\end{equation}
By using the Ritz projection $R_h$ defined in \eqref{eq:Ritz_proj}, we can alternatively write 
\[
u_h=u_h(q)=S_hq=R_h(Sq)=R_hu(q),\quad q\in Q.
\]
The discrete problem is
\begin{equation*}
\begin{aligned}
    \min_{q\in Q_{\mathrm{ad}},\,u_h\in V_h}\quad&
    J(q,u_h)
    =\frac12\|u_h-u_d\|_{L^2(\Omega)}^2
     +\frac{\alpha}{2}\|q\|_{L^2(\gamma)}^2,\\
    \text{subject to}\quad&
    (\nabla u_h,\nabla\chi)
    =\int_\gamma q\,\chi\,d \Hg
    \qquad\forall\chi\in V_h .
\end{aligned}
\end{equation*}
Notice that the control variable is not explicitly discretized. 
 
Let $j_h(q):=J(q,u_h(q))=J(q,S_hq)$ be the discrete reduced functional,
and the discrete variational problem can be formulated as 
\begin{equation}\label{eq: minimize discrete}
Minimize\quad j_h(q),\quad q\in Q_{ad}.
\end{equation}
\begin{theorem}\label{thm:existence reduced discrete}
There exists a unique solution $\bar{q}_h\in Q_{ad}$ to the minimization problem \eqref{eq: minimize discrete}.
\end{theorem}
\begin{proof}
The proof is identical to that of Theorem \ref{thm:existence reduced}, since $S_h:Q\to V_h$
 is linear and continuous and $\alpha>0$ makes $j_h$ strictly convex.
\end{proof}
The derivatives of $j_h(q)$, similarly to the continuous case, are  
\begin{subequations}
\begin{equation}\label{eq: jh prime}
j_h'(q)(\delta q)=(u_h-u_d,\delta u_h)+\alpha\langle q,\delta q \rangle_\gamma
\end{equation}
and
\begin{equation}\label{eq: jh double prime}
j_h''(q)(\delta q,\delta q)=\|\delta u_h\|^2_{L^2(\Om)}+\alpha\|\delta q\|^2_{L^2(\gamma)}
\end{equation}
\end{subequations}
where $u_h=S_hq$ and $\delta u_h=S_h\delta q$.
Similarly to the continuous case, $j_h'$ can be rewritten in terms of an adjoint variable $z_h=z_h(q)\in V_h$, which is the solution to the discrete adjoint equation 
\begin{equation}\label{eq: discrete adjoint}
(\na z_h,\na \chi)=({u}_h-{u}_d,\chi)\quad \forall \chi \in V_h.
\end{equation}
Thus, similarly to the continuous case, \eqref{eq: jh prime} can be rewritten as
\begin{equation}\label{eq: jh prime2}
\begin{aligned}
j_h'(q)(\delta q)&=(\na z_h,\na \delta u_h)+\alpha\langle q,\delta q \rangle_\gamma\\
&=\langle z_h, \delta q\rangle_\gamma+\alpha\langle q,\delta q \rangle_\gamma=\langle z_h+\alpha q,\delta q \rangle_\gamma.
\end{aligned}
\end{equation}
\begin{lemma}\label{lemma: characterization of qh}
A control $\bar{q}_h\in Q_{ad}$ is the solution to the problem \eqref{eq: minimize discrete} if and only if
\[
j_h'(\bar{q}_h)(\delta q-\bar{q}_h)\geq 0,\quad \forall \delta q\in Q_{ad}.
\]
\end{lemma}
\subsection{Discrete optimality conditions}

Using Lemma \ref{lemma: characterization of qh}, exactly as in the continuous case, the optimal solution $(\bar{q}_h,\bar{u}_h,\bar{z}_h)$ satisfies:
\begin{subequations}

the state equation: $\bar{u}_h\in V_h$ 
\begin{equation}\label{eq: discrete state2}
(\na \bar{u}_h,\na \chi) = \langle \bar{q}_h,\chi\rangle_\gamma\quad \forall\chi \in V_h,
\end{equation}

the adjoint equation: $\bar{z}_h\in V_h$ 
\begin{equation}\label{eq: discrete adjoint2}
(\na \bar{z}_h,\na \chi)=(\bar{u}_h-{u}_d,\chi),\quad \forall \chi \in V_h,
\end{equation}
and 
the optimality condition: $\bar{q}_h\in Q_{ad}$ 
\begin{equation}\label{eq: discrete grad2}
 \langle\alpha \bar{q}_h +\bar{z}_h\mid_\gamma,\delta q-\bar{q}_h\rangle_\gamma\geq 0,\quad \forall \delta q\in Q_{ad}.
\end{equation}
\end{subequations}
Equivalently,
\[
\bar q_h
=
P_{[q_a,q_b]}
\left(-\frac1\alpha\bar z_h|_\gamma\right)
\quad\text{a.e. on }\gamma.
\]

\begin{theorem}\label{thm:control-error}
	Let $\Omega\subset\mathbb R^3$ be a bounded convex
	polyhedral domain, let $\gamma\Subset\Omega$ be a compact
	embedded $C^{1,1}$ curve with two endpoints, and let
	$u_d\in L^\infty(\Omega)$.
	Let $\bar q\in Q_{\mathrm{ad}}$ be the optimal control
	of the continuous problem and let
	$\bar q_h\in Q_{\mathrm{ad}}$ be the solution of
	\eqref{eq: minimize discrete}.
	
	For continuous piecewise linear finite elements on a
	quasi-uniform, shape-regular family of tetrahedral meshes,
	with exact integration over $\gamma$, we have
	\[
	\alpha\|\bar q-\bar q_h\|_{L^2(\gamma)}
	\le
	Ch^2\ell_h^3
	\left(
	\|\bar q\|_{L^\infty(\gamma)}
	+\|u_d\|_{L^\infty(\Omega)}
	\right)
	\]
	for all sufficiently small $h$.
	The constant $C$ depends only on $\Omega$, the geometry
	of $\gamma$, its distance from $\partial\Omega$,
	and the mesh regularity constants.
	In particular, $C$ is independent of $h$ and $\alpha$.
\end{theorem}
\begin{proof}
	Since $j_h$ is quadratic and
	$j_h''(q)(\delta q,\delta q)\ge
	\alpha\|\delta q\|_{L^2(\gamma)}^2$, we have
	\[
	\begin{aligned}
		\alpha\|\bar q-\bar q_h\|_{L^2(\gamma)}^2
		&\le
		j_h'(\bar q)(\bar q-\bar q_h)
		-j_h'(\bar q_h)(\bar q-\bar q_h)\\
		&\le
		j_h'(\bar q)(\bar q-\bar q_h)
		-j'(\bar q)(\bar q-\bar q_h)\\
		&=
		\langle z_h(\bar q)-\bar z,\bar q-\bar q_h\rangle_\gamma\\
		&\le
		\|\bar z-z_h(\bar q)\|_{L^2(\gamma)}
		\|\bar q-\bar q_h\|_{L^2(\gamma)}.
	\end{aligned}
	\]
	Here we used the continuous and discrete optimality
	conditions in the second step.
	It therefore remains to estimate
	$\|\bar z-z_h(\bar q)\|_{L^2(\gamma)}$.
	
	Set $\widetilde z_h:=R_h\bar z$ and decompose
	\[
	\bar z-z_h(\bar q)
	=
	(\bar z-\widetilde z_h)
	+
	(\widetilde z_h-z_h(\bar q)).
	\]
	
	To estimate the first term, choose fixed open subdomains
	\[
	\gamma\Subset\Omega_0\Subset\Omega_1\Subset\Omega.
	\]
	The interior maximum-norm estimate
	\cite[Theorem~5.1]{AHSchatz_LBWahlbin_1977a},
	local interpolation, and the global $L^2$ estimate for
	the Ritz projection give
	\[
	\|\bar z-\widetilde z_h\|_{L^\infty(\Omega_0)}
	\le
	C\ell_h h^{2-3/p}\|\bar z\|_{W^{2,p}(\Omega_1)}
	+
	Ch^2\|\bar z\|_{H^2(\Omega)}
	\]
	for every $p\ge4$, with constants independent of $p$.
	By \cref{thm: regularity}, it follows that
	\[
	\begin{aligned}
		\|\bar z-\widetilde z_h\|_{L^2(\gamma)}
		&\le
		\mathcal H^1(\gamma)^{1/2}
		\|\bar z-\widetilde z_h\|_{L^\infty(\Omega_0)}\\
		&\le
		Cp^2\ell_h h^{2-3/p}
		\left(
		\|\bar q\|_{L^\infty(\gamma)}
		+\|u_d\|_{L^\infty(\Omega)}
		\right).
	\end{aligned}
	\]
	For sufficiently small $h$, we may choose $p=\ell_h$.
	Since
	\[
	h^{-3/p}
	=
	\exp\left(\frac{3|\ln h|}{\ell_h}\right)
	\le e^3,
	\]
	we obtain
	\[
	\|\bar z-\widetilde z_h\|_{L^2(\gamma)}
	\le
	Ch^2\ell_h^3
	\left(
	\|\bar q\|_{L^\infty(\gamma)}
	+\|u_d\|_{L^\infty(\Omega)}
	\right).
	\]
	
	For the second term, set
	$w_h:=\widetilde z_h-z_h(\bar q)$.
	The continuous and discrete adjoint equations imply
	\[
	(\nabla w_h,\nabla\chi)
	=
	(\bar u-u_h(\bar q),\chi)
	\qquad\forall\chi\in V_h.
	\]
	Testing with $\chi=w_h$, using the Cauchy--Schwarz
	inequality and \cref{lem:weighted-hardy}, yields
	\[
	\begin{aligned}
		\|\nabla w_h\|_{L^2(\Omega)}^2
		&=
		(\bar u-u_h(\bar q),w_h)\\
		&\le
		\|\sigma_h(\bar u-u_h(\bar q))\|_{L^2(\Omega)}
		\|\sigma_h^{-1}w_h\|_{L^2(\Omega)}\\
		&\le
		C\ell_h
		\|\sigma_h(\bar u-u_h(\bar q))\|_{L^2(\Omega)}
		\|\nabla w_h\|_{L^2(\Omega)}.
	\end{aligned}
	\]
	Consequently,
	\[
	\|\nabla w_h\|_{L^2(\Omega)}
	\le
	C\ell_h
	\|\sigma_h(\bar u-u_h(\bar q))\|_{L^2(\Omega)}.
	\]
	The discrete Sobolev inequality of 
	\cref{lem:discrete-sobolev} and
	\cref{thm: weighted error in L2} now give
	\[
	\begin{aligned}
		\|w_h\|_{L^2(\gamma)}
		&\le C\ell_h^{1/2}\|\nabla w_h\|_{L^2(\Omega)}\\
		&\le
		C\ell_h^{3/2}
		\|\sigma_h(\bar u-u_h(\bar q))\|_{L^2(\Omega)}\\
		&\le Ch^2\ell_h^3\|\bar q\|_{L^2(\gamma)}\\
		&\le Ch^2\ell_h^3\|\bar q\|_{L^\infty(\gamma)}.
	\end{aligned}
	\]
	Combining the estimates for the two terms, we conclude that
	\[
	\|\bar z-z_h(\bar q)\|_{L^2(\gamma)}
	\le
	Ch^2\ell_h^3
	\left(
	\|\bar q\|_{L^\infty(\gamma)}
	+\|u_d\|_{L^\infty(\Omega)}
	\right).
	\]
	Substitution into the initial inequality and division by
	$\|\bar q-\bar q_h\|_{L^2(\gamma)}$, unless it vanishes,
	prove the assertion.
\end{proof}
\begin{remark}\label{rem:piecewise-curves}
The regularity and error estimates extend to simple curves
$\gamma\Subset\Omega$ consisting of finitely many regular, embedded
$C^{1,1}$ arcs, including non-self-intersecting polygonal chains.
The weighted Hardy and discrete trace inequalities follow
by applying the corresponding estimates on each arc and summing.
The weighted state error estimate follows similarly by decomposing
the source into its arc contributions and using linearity.
Consequently, the control error estimate retains the rate
$h^2(1+|\ln h|)^3$. 
\end{remark}

\section{Numerical example for the variational discretization}
\label{sec:Numerics}

We consider
\[
  \Omega=(-1,1)^3,
  \qquad
  \gamma=\{(0,0,s):-\tfrac12\le s\le\tfrac12\},
\]
with regularization parameter and control bounds
\[
  \alpha=10^{-2},
  \qquad q_a=0,
  \qquad q_b=1.
\]
The desired state is the smooth function
\begin{equation}
  u_d(x,y,z)
  =5\exp\bigl(-25(x^2+y^2)-32z^2\bigr).
  \label{eq:numerical-desired-state}
\end{equation}

We use a sequence of nested, quasi-uniform tetrahedral meshes.  Each
coordinate direction is divided into $N$ subintervals, and every Cartesian
cube is subdivided into six tetrahedra sharing its body diagonal.  Thus,
\[
  h_N=\frac{2\sqrt{3}}{N}.
\]
For the values $N=4,8,\ldots,256$ used below, the segment $\gamma$ lies on
the mesh skeleton.  All curve integrals are evaluated on the exact curve;
hence no geometric approximation of $\gamma$ is introduced.

The control is not discretized independently.  Instead, the discrete
variational inequality gives the pointwise projection formula
\begin{equation}
  \bar q_h(s)
  =P_{[q_a,q_b]}
    \left(-\frac{\bar z_h(0,0,s)}{\alpha}\right),
  \qquad -\tfrac12\le s\le\tfrac12.
  \label{eq:numerical-projection}
\end{equation}
The discrete optimality system is solved by a primal--dual
active-set method. The iteration is terminated when both the relative
state--adjoint change and the relative change of the active-set operator
are below $10^{-8}$. We additionally require the relative block-system,
state, and adjoint residuals to be below $10^{-8}$.
 Each geometric edge of $\gamma$ is assembled exactly once. It is split
at the switching points of~\eqref{eq:numerical-projection}, so the active
and inactive contributions to the curve integrals are evaluated exactly.

Let $(\bar q_N,\bar u_N,\bar z_N)$ denote the discrete optimal
solution on the mesh with $N$ subdivisions in each coordinate direction.
Since an exact optimal control is unavailable, we compare controls on
successive nested meshes and define
\begin{equation}
  D_q(N):=
  \|\bar q_N-\bar q_{2N}\|_{L^2(\gamma)}.
  \label{eq:successive-control-difference}
\end{equation}
For the evaluation of $D_q(N)$, we use the common partition of
$\gamma$ containing the mesh points and the switching points of both
$\bar q_N$ and $\bar q_{2N}$. On each resulting interval their difference
is affine, and hence its squared $L^2$ integral is evaluated exactly.

The corresponding raw convergence rate is
\begin{equation}
  \operatorname{rate}_q(N)
  :=\log_2\!\left(\frac{D_q(N/2)}{D_q(N)}\right).
  \label{eq:successive-control-rate}
\end{equation}
Thus, $D_q(N)$ is a successive-mesh difference rather than an error against
the exact control.  Let
\[
 E_q(N):=\|\bar q-\bar q_N\|_{L^2(\gamma)}.
\]
Under the standard saturation assumption
$E_q(2N)\le\vartheta E_q(N)$ with some $\vartheta<1$, the triangle
inequality gives
\[
 (1-\vartheta)E_q(N)
 \le D_q(N)
 \le(1+\vartheta)E_q(N).
\]
Thus, in the asymptotic regime, $D_q(N)$ is equivalent to the unknown
control error up to mesh-independent constants.

\begin{table}[htbp]
  \centering
  \caption{Successive-mesh $L^2(\gamma)$ differences and raw convergence
  rates for the variationally discretized optimal control.}
  \label{tab:control-convergence}
  \begin{tabular}{ccc}
    \hline
    $N$--$2N$ & $D_q(N)$ & $\operatorname{rate}_q(N)$ \\
    \hline
    $4$--$8$     & $5.4913\times10^{-2}$ & --       \\
    $8$--$16$    & $6.7839\times10^{-2}$ & $-0.305$ \\
    $16$--$32$   & $2.7634\times10^{-2}$ & $1.296$  \\
    $32$--$64$   & $7.5384\times10^{-3}$ & $1.874$  \\
    $64$--$128$  & $1.9697\times10^{-3}$ & $1.936$  \\
    $128$--$256$ & $4.9859\times10^{-4}$ & $1.982$  \\
    \hline
  \end{tabular}
\end{table}

The first two mesh pairs exhibit preasymptotic behavior. Thereafter,
the raw rate increases toward two and reaches $1.982$ on the finest
mesh pair. For the fixed data of this experiment, Theorem~\ref{thm:control-error} gives
\[
  \|\bar q-\bar q_h\|_{L^2(\gamma)}
  \le Ch^2\ell_h^3,
  \qquad
  \ell_h=1+\abs{\ln h}.
\]

The observed successive-mesh rates are therefore consistent with the
almost quadratic estimate. Since neither an exact control nor an
asymptotic error expansion is available, this experiment supports the
predicted algebraic order but cannot identify or verify the power
$3$ of the logarithmic factor.

\section{Appendix: Proofs of the tubular neighborhood estimates}
\label{app:curve}
\begin{proof}[Proof of Lemma~\ref{lem:curve-tubular-coordinates}]
	Extend $\Gamma$ along its endpoint tangents by setting
	\[
	\widetilde\Gamma(s):=
	\begin{cases}
		\Gamma(0)+s\Gamma'(0), & -\delta\le s<0,\\
		\Gamma(s), & 0\le s\le L,\\
		\Gamma(L)+(s-L)\Gamma'(L), & L<s\le L+\delta.
	\end{cases}
	\]
	This extension is of class $C^{1,1}$ and is parametrized
	by arclength.
	By embeddedness and compactness of $\gamma$, and since
	$\gamma\Subset\Omega$, we can choose $\delta>0$ sufficiently
	small that $\widetilde\Gamma$ is injective and
	$\widetilde\gamma\Subset\Omega$.
	Fix an open interval
	\[
	[0,L]\subset J\Subset(-\delta,L+\delta).
	\]
	
	For the extended $C^{1,1}$ curve, one can directly construct
	Lipschitz continuous vector fields
	$n_1,n_2:\overline J\to\mathbb R^3$ such that
	$(\widetilde\Gamma'(s),n_1(s),n_2(s))$ forms a positively
	oriented orthonormal basis of $\mathbb R^3$ for every
	$s\in\overline J$.
	This also follows from the more general existence result
	for relatively parallel frames in
	\cite[Theorem~A.1]{KrejcirikZahradova2020}.
	In particular, $\Psi$ is Lipschitz continuous on
	$J\times B_R^{(2)}(0)$ for every fixed $R>0$.
	
	Using the orthogonality of
	$\Psi(s,y)-\widetilde\Gamma(s)$ to $\widetilde\Gamma'(s)$,
	the $C^{1,1}$ regularity, and embeddedness, one can
	directly verify that, for sufficiently small $R>0$,
	$\widetilde\Gamma(s)$ is the unique nearest point on
	$\widetilde\gamma$ to $\Psi(s,y)$ whenever
	$s\in J$ and $|y|<R$.
	Alternatively, this follows from the local projection
	theorem \cite[Theorem~1]{LeobacherSteinicke2021},
	applied to
	$M:=\widetilde\Gamma((-\delta,L+\delta))$,
	together with compactness of
	$\widetilde\Gamma(\overline J)$ and its positive distance
	from the endpoints of $\widetilde\gamma$.
	This proves \eqref{eq:tubular-distance}.
	Moreover, the inclusion $\gamma\subset\widetilde\gamma$
	immediately gives \eqref{eq:tubular-distance-comparison}.
	
	To see that $\Psi$ is injective, suppose
	$\Psi(s,y)=\Psi(t,z)$.
	Uniqueness of the nearest point implies
	$\widetilde\Gamma(s)=\widetilde\Gamma(t)$, and hence $s=t$.
	Orthonormality of $n_1(s),n_2(s)$ then gives $y=z$.
	By invariance of domain, the image $W^0$ is open.
	Since $\widetilde\Gamma(\overline J)\Subset\Omega$,
	we may decrease $R$ so that $W^0\Subset\Omega$.
	
	It remains to establish the Jacobian estimate.
	Set $T:=\widetilde\Gamma'$ and
	\[
	E(s,\vartheta)
	:=\cos\vartheta\,n_1(s)+\sin\vartheta\,n_2(s).
	\]
	Then, almost everywhere,
	\[
	\partial_s\Phi=T+\rho\,\partial_s E,
	\qquad
	\partial_\rho\Phi=E,
	\qquad
	\partial_\vartheta\Phi=\rho\,\partial_\vartheta E.
	\]
	Since $E\times\partial_\vartheta E=T$ and
	$\partial_s E\cdot T=-E\cdot T'$, it follows that
	\[
	|\det D\Phi|
	=\rho\,|1-\rho E\cdot T'|.
	\]
	Choosing $R$ still smaller if necessary so that
	\[
	R\|\widetilde\Gamma''\|_{L^\infty(J)}\le\frac12
	\]
	proves \eqref{eq:tubular-jacobian}.
	Furthermore, $|\partial_\rho\Phi|=|E|=1$.
	
	Finally, $W^0$ is an open neighborhood of the compact
	set $\gamma$.
	We can therefore choose an open set $W$ with
	$\gamma\subset W\Subset W^0$, and then $r>0$
	sufficiently small that
	\[
	\overline{T_{2r}(\gamma)}
	\subset W\Subset W^0\Subset\Omega.
	\]
	This proves \eqref{eq:tubular-cover}.
\end{proof}

\begin{proof}[Proof of Lemma~\ref{lem:curve-tubular-integration}]
	Choose $r_0>0$ smaller than $\operatorname{reach}(\gamma)$,
	the radius $R$ from Lemma~\ref{lem:curve-tubular-coordinates},
	and $\operatorname{dist}(\gamma,\partial\Omega)$.
	Then $T_{r_0}(\gamma)\Subset\Omega$, and every
	$x\in T_{r_0}(\gamma)$ has a unique nearest point
	$p(x)\in\gamma$.
	
	Fix $0<r\le r_0$ and set
	\[
	A_r:=\{x\in T_r(\gamma):
	p(x)\in\Gamma((0,L))\}.
	\]
	If $x\in A_r$ and $p(x)=\Gamma(s)$, the minimizing
	property implies
	\[
	(x-\Gamma(s))\cdot\Gamma'(s)=0.
	\]
	Thus, for $x\notin\gamma$, we can write
	$x=\Phi(s,\rho,\vartheta)$ with
	$s\in(0,L)$ and $\rho=d_\gamma(x)<r$.
	Conversely, the nearest-point property from
	Lemma~\ref{lem:curve-tubular-coordinates}, together with
	$\gamma\subset\widetilde\gamma$, gives
	\[
	p(\Phi(s,\rho,\vartheta))=\Gamma(s),
	\qquad
	d_\gamma(\Phi(s,\rho,\vartheta))=\rho
	\]
	for $s\in(0,L)$ and $0<\rho<r$.
	
	Consequently, after removing the curve itself and the
	surface corresponding to $\vartheta=0$, both of which
	have three-dimensional measure zero, $\Phi$ maps
	$(0,L)\times(0,r)\times(0,2\pi)$ bijectively onto
	the remaining part of $A_r$.
	The change-of-variables formula therefore yields
	\[
	\begin{aligned}
		\int_{A_r} F(d_\gamma(x))\,dx
		&=\int_0^L\int_0^{2\pi}\int_0^r
		F(d_\gamma(\Phi(s,\rho,\vartheta)))
		|\det D\Phi(s,\rho,\vartheta)|
		\,d\rho\,d\vartheta\,ds\\
		&=\int_0^L\int_0^{2\pi}\int_0^r
		F(\rho)\,
		|\det D\Phi(s,\rho,\vartheta)|
		\,d\rho\,d\vartheta\,ds\\
		&\le \frac32\int_0^L\int_0^{2\pi}\int_0^r
		F(\rho)\,\rho\,d\rho\,d\vartheta\,ds\\
		&=3\pi L\int_0^r F(\rho)\,\rho\,d\rho,
	\end{aligned}
	\]
	where the inequality follows from
	\eqref{eq:tubular-jacobian}.
	
	It remains to estimate the contributions associated
	with the endpoints.
	For $a\in\{\Gamma(0),\Gamma(L)\}$, set
	\[
	E_{r,a}:=\{x\in T_r(\gamma):p(x)=a\}.
	\]
	On this set, $d_\gamma(x)=|x-a|$, and
	$E_{r,a}\subset B_r(a)$.
	Since $F$ is nonnegative, spherical coordinates give
	\[
	\begin{aligned}
		\int_{E_{r,a}} F(d_\gamma(x))\,dx
		&\le \int_{B_r(a)} F(|x-a|)\,dx\\
		&=\int_{\mathbb S^2}\int_0^r
		F(\rho)\,\rho^2\,d\rho\,d\mathcal H^2(\omega)\\
		&=4\pi\int_0^r F(\rho)\,\rho^2\,d\rho\\
		&\le 4\pi r_0\int_0^r F(\rho)\,\rho\,d\rho.
	\end{aligned}
	\]
	The sets $A_r$, $E_{r,\Gamma(0)}$, and $E_{r,\Gamma(L)}$
	form a disjoint partition of $T_r(\gamma)$.
	Adding their contributions proves
	\eqref{eq:curve-tubular-integration}.
	Finally, taking $F\equiv1$ gives
	\eqref{eq:curve-tube-volume}.
\end{proof}


\clearpage
\bibliography{references}
\bibliographystyle{abbrvnat}

\end{document}